\documentclass[preprint,12pt]{elsarticle}

\usepackage{url}
\usepackage{graphicx}%
\usepackage{multirow}%
\usepackage{amsmath,amssymb,amsfonts}%
\usepackage{amsthm}%
\usepackage{mathrsfs}%
\usepackage[title]{appendix}%
\usepackage{xcolor}%
\usepackage{textcomp}%
\usepackage{manyfoot}%
\usepackage{booktabs}%
\usepackage{algorithm,algorithmic}%
\usepackage{listings}%

\usepackage{times}
\usepackage{subfigure}
\usepackage{caption}
\usepackage{amsmath,amsfonts,amssymb,mathtools,amsthm}
\usepackage{makecell}

\DeclareMathOperator*{\argmin}{arg\,min}
\usepackage{threeparttable}
\usepackage{threeparttablex}
\usepackage{minitoc} 

\usepackage{color, colortbl} 
\definecolor{lightgray}{gray}{0.9} 
\definecolor{lightgreen}{rgb}{0.8, 1.0, 0.8}  
\definecolor{lightblue}{rgb}{0.8, 0.9, 1.0}   
\definecolor{lightpink}{rgb}{1.0, 0.8, 0.9}  

\newcommand\red[1]{\textcolor{black}{#1}}
\definecolor{myblue}{RGB}{0,0,255}
\newcommand{\blue}[1]{\textcolor{myblue}{#1}}

\newtheorem*{prop*}{Proposition}

\newtheorem*{theorem*}{Theorem}

\newtheorem*{lemma*}{Lemma}

\newtheorem*{property*}{Property}

\newtheorem*{assumption*}{Assumption}

\usepackage{tabularx}
\usepackage{multirow}
\usepackage{booktabs}
\usepackage{enumitem}
\usepackage{lineno}

\journal{Computer Physics Communications}

\begin{document}

\begin{frontmatter}



	\title{Solving multiscale dynamical systems by deep learning}

	\author[1]{Junjie Yao\fnref{equal}} 
	\author[1]{Yuxiao Yi\fnref{equal}} 
	\author[1]{Liangkai Hang} 
	\author[3,4]{Weinan E} 
	\author[5,6]{Weizong Wang} 
	\author[1,2]{Yaoyu Zhang\corref{cor1}} 
	\ead{zhyy.sjtu@sjtu.edu.cn}
	\author[5,7]{Tianhan Zhang\corref{cor1}} 
	\ead{thzhang@buaa.edu.cn}
	\author[1,8]{Zhi-Qin John Xu\corref{cor1}} 
	\ead{xuzhiqin@sjtu.edu.cn}

	\affiliation[1]{organization={Institute of Natural Sciences, School of Mathematical Sciences,  MOE-LSC, Shanghai Jiao Tong University},city={Shanghai},postcode={200240},country={China}}
	\affiliation[2]{organization={Shanghai Center for Brain Science and Brain-Inspired Technology},city={Shanghai},postcode={201602},country={China}}
	\affiliation[3]{organization={Center for Machine Learning Research, School of Mathematical Sciences, Peking University},city={Beijing},postcode={100871},country={China}}

	\affiliation[4]{organization={AI for Science Institute},city={Beijing},country={China}}

	\affiliation[5]{organization={School of Astronautics, Beihang University},city={Beijing},postcode={100191},country={China}}
	\affiliation[6]{organization={State Key Laboratory of High-Efficiency Reusable Aerospace Transportation Technology}, city={Beijing}, postcode={102206}, country={China}}

	\affiliation[7]{organization={Key Laboratory of Spacecraft Design Optimization and Dynamic Simulation Technologies, Ministry of Education},postcode={102206}, city={Beijing}, country={China}} 

	\affiliation[8]{organization={Key Laboratory of Marine Intelligent Equipment and System, Ministry of Education},postcode={201100},country={China}} 

	\fntext[equal]{These authors contributed equally to this work and listed in alphabetical order.}
	\cortext[cor1]{Corresponding author.}
	\begin{abstract}
		Multiscale dynamical systems, modeled by high-dimensional stiff ordinary differential equations (ODEs) with wide-ranging characteristic timescales, arise across diverse fields of science and engineering, but their numerical solvers often encounter severe efficiency bottlenecks.
		This paper introduces a novel \textbf{DeePODE} method, which consists of an Evolutionary Monte Carlo Sampling method (EMCS) and an efficient end-to-end deep neural network (DNN) to predict multiscale dynamical systems. The method's primary contribution is its approach to the ``curse of dimensionality''-- the exponential increase in data requirements as dimensions increase. By integrating Monte Carlo sampling with the system's inherent evolutionary dynamics, DeePODE efficiently generates high-dimensional time-series data covering trajectories with wide characteristic timescales or frequency spectra in the phase space. Appropriate coverage on the frequency spectrum of the training data proves critical for data-driven time-series prediction ability, as neural networks exhibit an intrinsic learning pattern of progressively capturing features from low to high frequencies. We validate this finding across dynamical systems from ecological systems to reactive flows, including a predator-prey model, a power system oscillation, a battery electrolyte thermal runaway, and turbulent reaction-diffusion systems with complex chemical kinetics. The method demonstrates robust generalization capabilities, allowing pre-trained DNN models to accurately predict the behavior in previously unseen scenarios, largely due to the delicately constructed dataset. While theoretical guarantees remain to be established, empirical evidence shows that DeePODE achieves the accuracy of implicit numerical schemes while maintaining the computational efficiency of explicit schemes. This work underscores the crucial relationship between training data distribution and neural network generalization performance. This work demonstrates the potential of deep learning approaches in modeling complex dynamical systems across scientific and engineering domains.
	\end{abstract}


\end{frontmatter}


\doparttoc 
\faketableofcontents 
\part{} 
\vspace{-3em} 




\section{Introduction}\label{sec1}

Multiscale dynamical systems frequently arise in various scientific problems, including chemical kinetics \cite{barlow1998effects}, generic circuit simulation \cite{mcadams1995circuit}, ecosystem evolution \cite{volterra1931variations}, oscillation in power systems \cite{ye2016analysis}, and biological neural network modeling \cite{sun2009library}. Such a multiscale dynamical system often exhibits widespread characteristic timescales and can be modeled by high-dimensional stiff ordinary differential equations (ODEs). Classical numerical solvers require a small time step to guarantee numerical stability, resulting in high computational costs.
Consequently, over the decades, a sustained effort has been made to develop efficient and robust numerical methods to simulate multiscale systems across diverse scientific disciplines.

Machine learning approaches have been introduced to solve high-dimensional multiscale ODEs, but the neural network structure or limited training data often constrains their generality. 
For example, physics-informed neural networks \cite{dissanayake1994neural,raissi2019physics}  require a large training cost but only apply for a specific initial condition, preventing their applications in complex spatiotemporal systems such as turbulent flames~\cite{zhang2022multi} and biological neural networks~\cite{li2023mathematical}. 
Another frequently employed machine learning methodology focuses on the state mapping function~\cite{LED2022, AdaLED2023, GLED2024}, whose input is the current state and output is the subsequent state or the state change after a specific small or even large time step size. This approach can enhance generality once sufficient training data covers the dynamical evolution trajectories in phase space.

Capturing representative patterns in high-dimensional multiscale dynamical systems presents a fundamental challenge in data-driven modeling, especially for enabling the surrogate model to generalize effectively across diverse scenarios.
The difficulty stems from two key factors:  \textit{1) the "curse of dimensionality"} where data requirements grow exponentially with system dimensions, and \textit{ 2) the presence of multiple timescales} which demands sampling across vastly different temporal resolutions. While traditional Monte Carlo (MC) methods can theoretically overcome the dimensionality challenge, they do so at the cost of extremely low efficiency. In the case of MC sampling, the error decreases with the sample number $M$ at a rate $1/\sqrt{M}$. As an example for illustration, the concentration of measure shows that uniformly random samples concentrate on the surface of the high-dimensional sphere \cite{milman2009asymptotic}. Tabulation approach \cite{chen1989pdf,pope1997computationally} or manifold sampling (sampling from a set of evolution trajectories) in reactive flow simulation \cite{christo1996artificial,blasco1998modelling,choi2005fast,sen2010linear,sinaei2017large,zhang2021deep,owoyele2022chemnode,almeldein2022accelerating,de2022physics,yao2022gradient} mainly tackle such problem by dimension reduction, which can only work for low-dimensional systems or few working conditions. Successful examples often utilize the intrinsic characteristics of considered problems. For example, the key in AlphaGo for solving the Go problem is the high-efficiency sampling in high-dimensional space by MC tree search \cite{silver2016mastering}.

Unlike intractable complexity of many problems, such as the Bellman equation in AlphaGo, multiscale systems governed by explicit ODEs offer a significant advantage (\textit{local dynamical behavior}): their governing equations provide valuable information about the target function's behavior across high-dimensional space. This inherent mathematical structure can be leveraged to develop more efficient solutions. 

In this work, we propose a novel DeePODE method to establish global surrogate model for efficiently predicting multiscale dynamical systems. The DeePODE method takes on the \textit{"local dynamical behavior"} advantage through two components: the evolutionary Monte Carlo sampling method (EMCS) and the end-to-end deep neural network (DNN). The EMCS method combines the strengths of traditional Monte Carlo approaches with the system's natural evolutionary dynamics. More specifically, the method uses Monte Carlo sampling to achieve global coverage of the high-dimensional space, addressing the curse of dimensionality. This global sampling is then enhanced by evolving each sample point along its local ODE trajectory at a speed determined by the system's characteristic timescales. 
This hybrid approach naturally incorporates multiscale information into the sampling process. In regions where the system evolves slowly (low gradient regions), the method automatically collects fewer samples since large time steps are adopted to depict such domains. This hybrid sampling strategy efficiently generates representative data across the global high-dimensional space while respecting the system's inherent behavior. Then, the DNN is designed to predict system dynamics using arbitrarily reasonable large time steps. This network embeds temporal evolution patterns of wide frequency spectrum, enabling accurate predictions over extended time intervals that would be computationally expensive with conventional numerical methods. This work focuses on demonstrating the importance of sampling representative data, which enable neural networks, even a vanilla fully-connected network, to achieve high accuracy in simulating high-dimensional and stiff ODEs.

We validate the effectiveness of DeePODE in various multiscale systems, including a predator-prey model, an electronic dynamical process, a battery electrolyte auto-ignition, and several turbulent reaction-diffusion systems considering detailed chemical kinetics. These systems represent critical areas of scientific and engineering interest with wide-ranging applications, in which DeePODE achieves satisfactory accuracy, robustness and computational efficiency.
One particularly challenging area where DeePODE excels is in reaction-diffusion system (\textit{i.e.} combustion) simulations \cite{hawkes2005direct}. The complexity of these simulations stems from multiple factors: the high dimensionality of chemical reaction processes involving numerous species, time scales spanning eight orders of magnitude, and the influence of turbulent flow fields. Such system is crucial for understanding coupled reaction and transport phenomenon, from propulsion system to atmospheric chemistry and interstellar reactions. 

DeePODE demonstrates remarkable capabilities in addressing these challenges. One of its key strengths lies in the global surrogate it provides for combustion simulation—a single pretrained ODE surrogate model that can seamlessly couple with 1D, 2D, or 3D PDEs under previously unseen initial conditions, without requiring retraining or fine-tuning. Another advantage is the adaptability of the model. DeePODE models can be readily integrated into different numerical platforms, offering significant computational speed improvements—often two orders of magnitude faster than traditional solvers. This combination of accuracy, acceleration, and adaptability makes DeePODE a powerful tool for complex system simulations, such as biomass conversions \cite{mettler2012top} and air pollution predictions \cite{seinfeld1998air}.

\section{Methodology}\label{sec2}

In this section, we introduce DeePODE by demonstrating its two parts: evolutionary Monte Carlo sampling and deep neural network prediction.

\begin{figure*}[htbp]
\centering
    \includegraphics[width=1\textwidth]{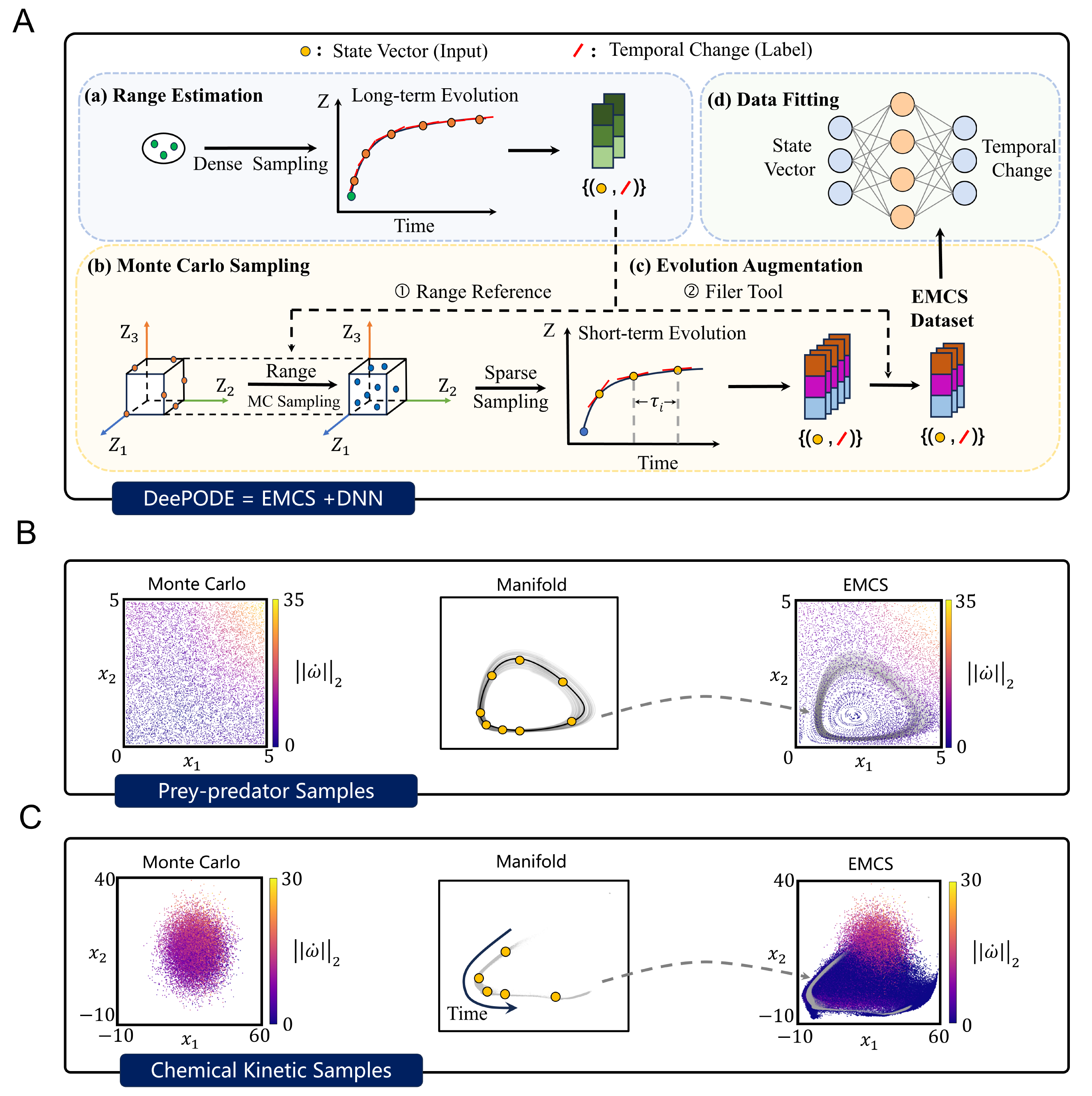} 
\caption{\textbf{Schematic diagram of the DeePODE method}. \textbf{(A)} Flowchart of DeePODE method. After obtaining a manifold dataset and an MC dataset, each point of the MC dataset is taken as the initial condition of a corresponding dynamical system. Sample along the evolution trajectory of MC data sparsely and filter them by the temporal gradient range of the manifold dataset. The remaining data composes the final EMCS dataset. A DNN is used to fit the dataset. Our DeePODE method can be regarded as EMCS combined with DNN, namely "DeePODE = EMCS + DNN" as in the annotation.
\textbf{(B)} Example of prey-predator system. The left and right figures represent the data distribution obtained by the Monte-Carlo and the EMCS methods, respectively, and the color of the data points represents the change's norm. The smooth gray regions in the second and third figures represent some evolution trajectories of this system. The yellow circular markers denote the sampled points in the EMCS method. \textbf{(C)} Example of Chemical kinetic system. Three figures serve the same purpose as the corresponding figures in \textbf{(B)} except that $x_1$ and $x_2$ are the data projected onto the first two principal directions of the manifold dataset. (For interpretation of the colors in the figure(s), the reader is referred to the web version of this article.) } \label{fig:sampling}
\end{figure*}

\subsection{Problem description}
Consider the high-dimensional dynamical system of the form
\begin{align}
    \frac{\mathrm{d}\boldsymbol{x}}{\mathrm{d}t} = \boldsymbol{f}(\boldsymbol{x}, t)
\end{align}
The system's evolution process can be expressed as: 
\begin{align}
    \boldsymbol{x}(t+\Delta t) = \mathcal{F}(\boldsymbol{x}(t))
\end{align}
where $\boldsymbol{x}(t)\in\mathbb{R}^{d}$ is the d-dimensional state vector, $\Delta t$ is the desired time step and $\mathcal{F}: \mathbb{R}^d\rightarrow\mathbb{R}^{d}$ denotes the propagator function mapping the system state from one time point to the next. Traditional numerical solvers face a significant limitation: they must use very small time step $\delta t\ll\Delta t $ to maintain stability when dealing with stiff systems. This requirement leads to high computational cost. DeePODE method addresses this limitation by constructing a parameterized data-driven surrogate model $\mathcal{F}_{\theta}: \mathbb{R}^d\rightarrow\mathbb{R}^{d}$. This surrogate model can directly predict the system state over larger time steps, effectively replacing the computationally intensive step-by-step integration process with a more efficient end-to-end single-step prediction. 

\subsection{Evolutionary Monte Carlo sampling}

Developing an end-to-end surrogate model is conceptually straightforward, but its success heavily depends on the quality and representativeness of training data. Our key innovation lies in the data generation technique, specifically the Evolutionary Monte Carlo Sampling (EMCS) method, which achieves comprehensive coverage through three sequential steps.

\paragraph{Step 1: Range estimation}

Due to the unique nature of each ODE system, it is important to understand its characteristic data distribution at the beginning. The range estimation process addresses this need by determining the system's feasible ranges through long-term evolution from several extreme initial conditions of interest. \red{To account for the diverse scales of dynamical systems (e.g., species concentrations in chemical kinetics), we randomly select physically realistic initial conditions based on the system’s physical constraints—for example, in chemical reaction kinetics, only fuel and oxidizer are present initially, with other species’ mass fractions set to zero, and temperatures are capped at realistic values (e.g., 3000 K). These conditions are evolved numerically using CVODE with a small time step (e.g., $10^{-8}$ s) until stabilization or periodic behavior emerges. The resulting evolutionary trajectories form the dataset $\mathbb{D}_{\text{init}} = \{ \langle \boldsymbol{x}(t), \boldsymbol{u}(t) \rangle \}$, where $\boldsymbol{u}(t) := \boldsymbol{x}(t + \Delta t) - \boldsymbol{x}(t)$ is the temporal change. This process, requiring only ~1\% of the full training data volume, provides an efficient and robust range estimation for subsequent sampling steps.
    As shown in Fig.~\ref{fig:sampling}A(a),} the dataset structure can be visualized through paired yellow dots and green columns. Each yellow dot represents a phase state $\boldsymbol{x}(t)$ in the high-dimensional state space, while its corresponding green column represents the temporal change \red{$\boldsymbol{u}(t):=\boldsymbol{x}(t+\Delta t)-\boldsymbol{x}(t)$, where $\Delta t$ is the time step used for DNN prediction}. Together, they form a training pair $(\boldsymbol{x}(t), \boldsymbol{u}(t))$, where $\boldsymbol{x}(t)$ serves as the input data point and $\boldsymbol{u}(t)$ as its label. By tracking these trajectories and collecting samples along them, we obtain a rough but essential estimation of the feasible range $[x^{\min}_{i}, x^{\max}_{i}]$ and $[u^{\min}_{i}, u^{\max}_{i}]$ for each component $x_i$ in the state vector $\boldsymbol{x}_t$. This initial dataset, denoted as $\mathbb{D}_{\text{Init}}:={(\boldsymbol{x}(t), \boldsymbol{u}(t))}$, provides crucial information about the system's bounds. A typical complex chemical reaction system like electrolyte thermal runaway would involve simulating electrolytes from various extreme temperature and concentration combinations to understand the full range of possible data distribution. This comprehensive range estimation is essential because electrolyte thermal runaway involves complex interactions between temperature-dependent decomposition reactions, heat generation, and species evolution, each operating across different scales and potentially leading to critical transitions in system behavior.

\paragraph{Step 2: Monte Carlo (MC) sampling}
After establishing the rough ranges in Step 1, we obtain a high-dimensional hypercube that bounds our system's state space. The Monte Carlo sampling then operates within this hypercube, with extra consideration for the data's scale distribution as shown in Fig. \ref{fig:sampling}A(b). When data components share similar orders of magnitude, we employ linear-scale sampling, such as the prey-predator model. However, for components spanning multiple orders of magnitude, we switch to logarithmic-scale sampling to ensure appropriate coverage across all scales, such as chemical reaction systems where species concentrations and reaction rates often vary by several orders of magnitude.

The dataset $\mathbb{D}_{\text{Init}}$ not only defines the sampling bounds but also guides the filtering of sampled points. We filter the Monte Carlo samples based on their rate of change components $\boldsymbol{u}(t)$, comparing them against the reference ranges observed in $\mathbb{D}_{\text{Init}}$. This filtering step is crucial for eliminating unphysical combinations that might arise from independent sampling of each dimension. For example, in chemical systems, certain species concentrations might be mathematically possible within the hypercube but physically impossible due to reaction constraints or conservation laws. The filter thus removes samples with unrealistic rates of change, ensuring the dataset maintains physical relevance.

\paragraph{Step 3: Evolution augmented generation}
The final step, evolution augmented generation, enhances the dataset by evolving each Monte Carlo sample point along its natural trajectory according to the ODEs as shown in Fig. \ref{fig:sampling}A(c). By sampling data points from these local trajectories, we effectively capture the system's behavior in the manifold surrounding each Monte Carlo point. A critical insight in this evolution step is the use of non-uniform evolution time steps  $\Delta t_i=\boldsymbol{\tau_i} = [\tau_1, \tau_2, \cdots, \tau_k]$, which enables the updated dataset $\mathbb{D}_{\rm EMCS}$ to capture more characteristic timescales of the system. By allowing the sampling intervals to vary based on local dynamics, we can efficiently represent both rapid changes and slow evolution within the system. \red{In practice, we employ an increasing $\tau_i$ strategy where $\tau_{i+1} \geq \tau_i$. The selection criteria for $\tau_i$ are discussed in detail in Section \ref{sec:csp_time_gms}.}

This approach simultaneously achieves two crucial objectives: broad phase space coverage for robust training and detailed local dynamics for accurate prediction. The combination of global Monte Carlo sampling with locally adapted trajectory evolution ensures that our dataset captures both the overall structure of the solution space and the fine-grained temporal dynamics at different scales. This comprehensive sampling strategy forms the foundation for building reliable surrogate models of complex dynamical systems.

\subsection{Deep neural network}

The feed-forward network (FFN), denoted as $\mathcal{F}_{\boldsymbol{\theta}}$, is leveraged to fit the training dataset, and the net structure is shown in Fig. \ref{fig:sampling}A(d). In the first hidden layer, the input $\boldsymbol{x}$ undergoes an affine transformation with trainable weights $\boldsymbol{W}^{\scriptscriptstyle{[1]}}$ and biases $\boldsymbol{b}^{\scriptscriptstyle{[1]}}$, followed by a GELU activation function, \textit{i.e.}, $\sigma(\boldsymbol{W}^{\scriptscriptstyle{[1]}} \boldsymbol{x} +\boldsymbol{b}^{\scriptscriptstyle{[1]}})$. Iteratively, the output of the previous hidden layer is used as the input of the subsequent hidden layer until the output layer. 

Generally, the input of FFN is the state vector $\boldsymbol{x}(t)$, and the output is the temporal change $\boldsymbol{u}(t)=\boldsymbol{x}(t+\Delta t)-\boldsymbol{x}(t)$, where $\Delta t=10^{-6}$ s unless otherwise specified. Specifically, data pre-processing is essential for combustion examples. We apply Box-Cox transformation (BCT) \cite{Box1964} to the mass fraction in $\boldsymbol{x}(t)$ before it is fed into the FFN, transforming it into $\mathcal{O}(1)$ order. Note that the BCT of $x$ is $(x^{\lambda}-1)/\lambda$ where $\lambda$ is set as 0.1. As indicated by the frequency principle  \cite{xu_training_2018,xu2019frequency}, the neural network is difficult to learn small-scale components, and BCT is an effective way to alleviate such difficulty~\cite{zhang2022multi}. The optimal parameters are determined by minimizing the mean absolute error: \begin{align}
    \boldsymbol{\theta}^{*} = \argmin_{\boldsymbol{\theta}}\frac{1}{N}\sum^{N}_{i=1}\vert \boldsymbol{u}_i - \mathcal{F}_{\boldsymbol{\theta}}(\boldsymbol{x}_i) \vert 
\end{align}
where $N$ denotes the size of the training dataset. The optimization algorithm is Adam with batch size 1024 and learning rate $10^{{-4}}$. 

\subsection{Sampling distribution comparison}
To illustrate EMCS's effectiveness, we compare data distributions generated by crude MC sampling and EMCS methods for two distinct systems: a prey-predator model and a complex chemical reaction system.

Fig. \ref{fig:sampling}(B) presents three views of the prey-predator system's data distribution. The left panel shows the uniform distribution from basic Monte Carlo sampling, which treats prey and predator populations as independent dimensions. The middle panel displays a characteristic oscillating trajectory obtained through long-term system evolution, forming a closed loop representing the system's natural attractor. The right panel reveals a key advantage of EMCS: it automatically captures the data pattern from the ODE system evolution. We will show in the Result part that the neural network trained on EMCS data can accurately predict the system's future states, while the MC-trained network fails due to the lack of representative samples.

For the complex chemical reaction system of methane oxidation, as shown in Fig. \ref{fig:sampling}(C), we use the first two principal components after principal component analysis (PCA) to visualize the complex system with more than 50 dimensions. The distribution of the PCA-projected MC samples is homogeneous in different directions, as shown in the left panel.

 The middle panel shows a typical smooth trajectory after PCA projection, considered a test manifold dataset. With the EMCS method evolving the governing equations, the obtained dataset, as shown on the right side, has more samples covering the region of the manifold dataset (gray region) than MC sampling.
Consequently, EMCS achieves comprehensive coverage across the high-dimensional phase space.

\section{Results}\label{sec3}

We test out the method on several tasks of increasing complexity, including 1) a two-dimensional autonomous predator-prey model, 2) a 15-dimensional non-autonomous electronic circuit system, 3) a 104-dimensional battery electrolyte thermal runaway, and 4) a series of challenging turbulent diffusion-reaction system simulations.

\subsection{Two-dimensional predator-prey model}\label{sec:predator-prey}

Consider a two-dimensional Lotka-Volterra predator-prey model \cite{qin2021data}: $\mathrm{d} x_1/\mathrm{d} t  =-x_1 + x_1 x_2$ and $\mathrm{d} x_2/\mathrm{d} t=2x_2 - x_1 x_2$.
We demonstrate DeePODE's sampling advantages using this classical example. The model architecture employs a deep neural network with three hidden layers, each containing 200 neurons. This network maps the system's current state to its future state after a time step of $\Delta t=0.1$ seconds. This simple example showcases how DeePODE effectively combines the benefits of both Monte-Carlo sampling (broad state space coverage) and evolution augmentation (concentration in dynamically important regions), achieving comprehensive coverage while maintaining accuracy in critical regions.

We compare three sampling strategies, each generating 100,000 data points. 1) Monte-Carlo (MC) sampling: We randomly sample 100,000 initial points within the domain $x_1\in[0,5]$ and $x_2\in[0,5]$. Due to system stiffness, labels are computed using direct integration with time steps under $10^{-2}$ seconds. 2) EMCS sampling: Starting with 20,000 random initial points, we augment the dataset by evolving each point through time steps $\boldsymbol{\tau} = [0.1, 0.2, 0.3, 0.4]$ seconds, yielding the full 100,000-point dataset. 3) Manifold sampling: We select 500 random initial points and sample 200 points along each trajectory at fixed intervals of 0.1 seconds, producing 100,000 points from 500 evolutionary trajectories. This sampling comparison maintains consistent dataset sizes while highlighting each method's distinct approach to state space exploration.

\begin{figure}[htbp]
	\centering
        \includegraphics[width=1\textwidth]{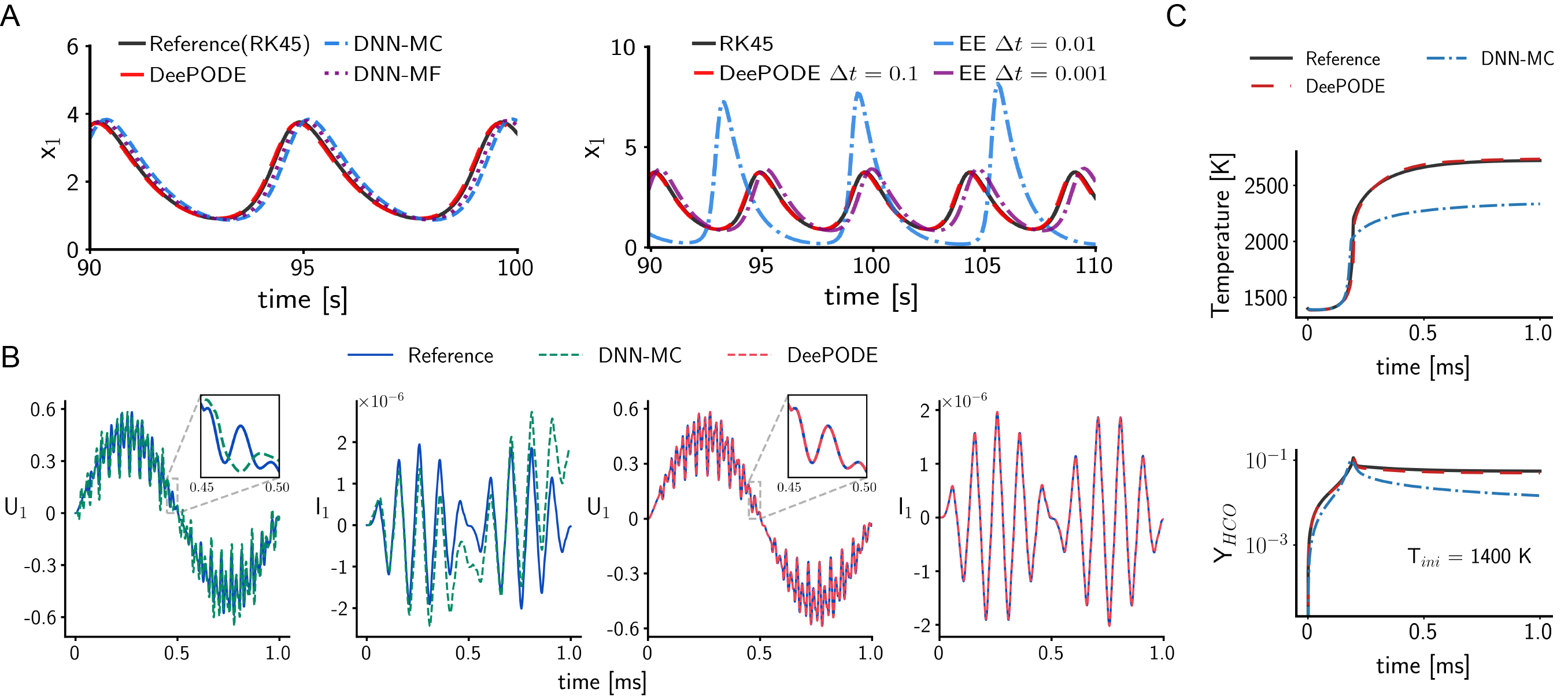}
\caption{\textbf{DeePODE for multiscale dynamic systems}.  \textbf{(A)} Two-dimensional predator-prey model. DNN prediction and direct integration results for $x_1$ from the initial value $x_1=3$, $x_2=2$. RK45, MC, MF, and EE represent 
Runge-Kutta-Fehlberg scheme, Monte Carlo, manifold sampling methods, and Explicit Euler scheme, respectively. \textbf{(B)} Ring modulator model. DNN prediction and direct integration results for $U_1$ and $I_2$. \textbf{(C)} Battery electrolyte temporal evolution of temperature and radical $\text{HCO}$ mass fraction. The initial condition is $p_0=1$ atm, $T_0=1400$ K, equivalence ratio $\phi_0=1$.
\label{fig:predator_modulator_battery} }

\end{figure}

 The performance comparison reveals distinct outcomes across different methods. As shown in Fig. \ref{fig:predator_modulator_battery}(A), the DNNs trained using the MC sampling or manifold sampling show substantial deviations from the reference solution. In contrast, DeePODE maintains accurate predictions even over extended time periods (up to t=100 seconds). \red{Although all three sampling methods use datasets of identical size, they produce distinct data distributions (Fig.~\ref{fig:sampling}(B)). This demonstrates that the spatial distribution of training data critically influences DNN performance, with EMCS achieving broader and more representative coverage of the phase space.}  Moreover, DeePODE achieves this accuracy while using larger integration steps than the explicit Euler method, demonstrating a crucial advantage: it combines the computational efficiency of explicit schemes with enhanced numerical stability. This combination of accuracy and efficiency makes DeePODE particularly valuable for long-time integration of dynamical systems.

\subsection{Electronic circuit system}
The ring modulator represents a challenging test case for numerical methods in electronic circuit simulation. This non-autonomous system incorporates 15 electronic components with highly nonlinear interactions and multiple time scales. Traditional numerical solvers require extremely small integration time steps (approximately $10^{-8}$ seconds) to maintain stability and accuracy due to the system's stiffness. This stringent time-step requirement makes conventional simulation approaches computationally expensive, particularly for long-time predictions. The complexity of this system, characterized by its high nonlinearity and multiscale nature, makes it an ideal benchmark for evaluating numerical method performance.

Our sampling strategy for the ring modulator system consists of two approaches for comparison. In the DeePODE approach, we initially sample 80,000 points within the system's natural bounds, determined from traditional solver results. We then evolve each initial point for $10^{-3}$ seconds, taking 30 random samples along each trajectory, resulting in a total dataset of 2,400,000 points. For comparison, we implement a purely random sampling approach by directly sampling 2,400,000 points within the same range. Both approaches use identical neural network architectures, consisting of three hidden layers with 800, 400, and 200 neurons, respectively. The network maps the input vector containing the current state and time to the system state after $\Delta t = 10^{-6}$ seconds. This setup allows direct comparison between evolution-augmented and pure random sampling while maintaining consistent dataset sizes and model architectures.

The simulation results in Fig. \ref{fig:predator_modulator_battery} (B) compare the performance of different methods in predicting two high-frequency signals: voltage $U_1$ and current $I_1$. Using the Runge-Kutta solution as our reference, DeePODE demonstrates remarkable accuracy in capturing high-frequency oscillations. In contrast, the DNN trained with Monte Carlo sampling (DNN-MC) shows significant deviations from the reference solution, despite using the same dataset size. This comparison underscores DeePODE's capability to handle high-dimensional, multiscale systems effectively, while pure Monte Carlo sampling struggles to achieve comparable accuracy. The superior performance of DeePODE highlights the importance of evolution-augmented sampling in capturing the system's dynamic behavior. Appendix \ref{sec:modulator_appd} shows more information about the case.

\subsection{Battery electrolyte thermal runaway}

In this section, we apply DeePODE to study the thermal runaway process in lithium batteries, focusing on the thermal decomposition of dimethyl carbonate (DMC) electrolyte. During thermal runaway, a critical battery safety concern, the electrolyte decomposes and generates flammable gases. The chemical system is modeled using a comprehensive mechanism containing 102 species and 805 elementary reactions, specifically developed to represent the gas mixture released during battery thermal runaway. We aim to develop a DNN-based surrogate model that accurately captures the complex thermal decomposition dynamics under constant pressure conditions. This application demonstrates DeePODE's ability to handle the intricate chemistry involved in battery safety analysis.

We utilize Cantera to obtain the reference solution and generate datasets for DNN models. More specifically, for the MC method, 500,000 data points are sampled uniformly within the range $\text{temperature}\in[300, 3000]$ K, $\text{pressure}\in[0.5,3]$ atm, $\text{species mass fraction} \in[0, 1]$.  For the DeePODE method, we randomly pick 100,000 data points within the same range 
and apply EMCS for those points with evolution time sequence $\boldsymbol{\tau}=[10^{-6}, 10^{-5}, 10^{-4}, 10^{-3}]$ s. The corresponding label datasets are obtained by evolving the input vectors for $\Delta t = 10^{-6}$ s. We use the two datasets to train the corresponding DNNs with 1600,800,400 hidden neurons for three hidden layers, respectively. With the stiffness constraint, the integration time step inside Cantera is $10^{-8}$, while the DNN models predict the state after $\Delta t=10^{-6}$ given an input state. 

The simulation results shown in Fig. \ref{fig:predator_modulator_battery}(C) demonstrate DeePODE's performance on a complex 104-dimensional chemical system, tracking both key species concentrations such as HCO radical and temperature evolution. Starting from initial conditions of 1400 K temperature, 1 atm pressure, and stoichiometric mixture ($\phi = 1.0$), The results show that DeePODE maintains high accuracy throughout the simulation, particularly during the critical thermal runaway period where temperature changes rapidly. In contrast, the Monte Carlo-trained model fails to accurately capture the system behavior during this crucial turning point, where the chemistry is most complex and the temperature gradient is steepest.

\subsection{Reaction-diffusion system with complex chemical kinetics}\label{sec:chemical}

\color{black}{We evaluate DeePODE using a complex reaction-diffusion system that combines chemical reactions with fluid dynamics. This system is a rigorous test case because it involves interactions across multiple scales and bidirectional coupling between processes. The chemical reactions generate heat and transform chemical species, while fluid motion redistributes both heat and species throughout the system. These processes create a challenging feedback loop: reaction-generated heat alters fluid movement, which then influences the transport of chemical species and subsequent reactions. The system's sensitivity to small perturbations makes it particularly suitable for testing numerical methods. Even minor inaccuracies in calculating reaction rates can compound over time, potentially leading to significant errors in predicting flame propagation. DeePODE demonstrates remarkable stability under these demanding conditions, maintaining accurate predictions despite the system's inherent complexity. This performance in handling tightly coupled, multiscale phenomena suggests that DeePODE offers a reliable solution for modeling complex physical systems.} 

In this section, we evaluate DeePODE models using several representative chemical kinetic models, including methane oxidation model GRI-Mech 3.0, which comprises 53 species and 325 reactions; DRM19, with 21 species and 84 reactions for methane/air reaction; and an n-heptane oxidation model containing 34 species and 191 reactions \cite{DeePMR2022}. We train a corresponding DeePODE model for each mechanism, which we have made publicly available through our repository (https://github.com/intelligent-algorithm-team/intelligent-combustion). These models are then thoroughly evaluated using various flow benchmark tests. A key strength of DeePODE is its versatility - once trained, a single model can be directly integrated into one-, two-, or three-dimensional flow simulations without requiring additional training or adjustment. This "train once, use anywhere" capability significantly enhances its practical reactor applications. 
The CFD code used for algorithm tests is \textit{EBI-DNS} \cite{zirwes2020quasi}. To demonstrate the portability of DeePODE, we integrate the DNN model into two widely-used CFD frameworks: \textit{ASURF} \cite{chen2009effects,zhang2021studies}  and \textit{DeepFlame} \cite{mao2022deepflame} codes. \red{The ODE solver used for both generating training data and evaluation is CVODE\cite{cohen1996cvode}, which employs the implicit Backward Differentiation Formula (BDF) method, designed for robust handling of stiff systems.} Further details on chemical kinetics sampling and additional examples are provided in the Appendix \ref{sec:combustion_appd}.

\begin{figure*}[htp]
	\centering
        \includegraphics[width=1\textwidth]{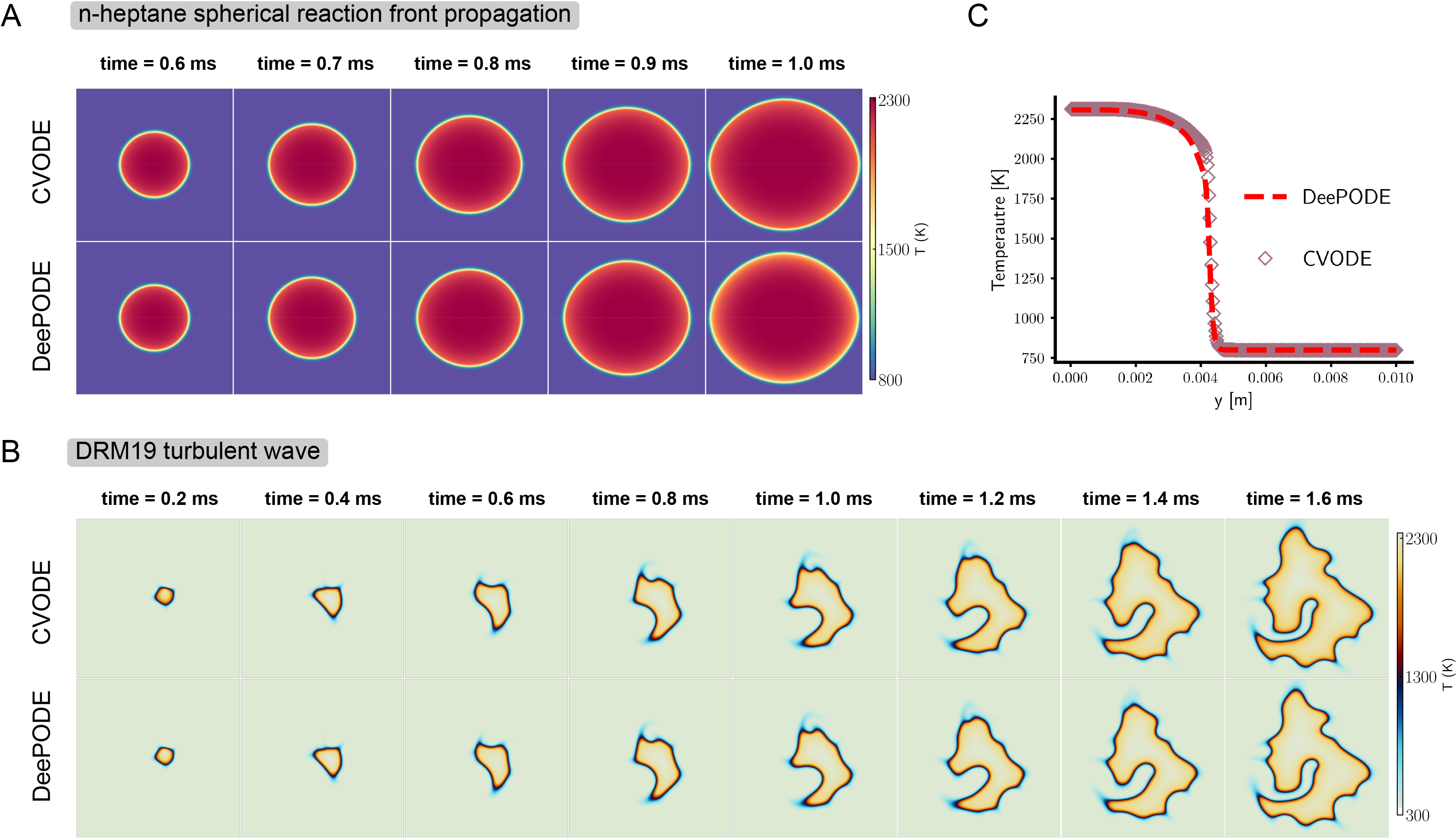}
    \caption{\textbf{Two-dimensional spherical case and turbulent case}. The snapshots of simulation comparison at different time instants are obtained with CVODE and DeePODE of \textbf{(A)} n-heptane spherical reaction wave and \textbf{(B)} DRM19 turbulent reaction-diffusion wave, respectively. Comparison of cross-section temperature profiles of \textbf{(C)} n-heptane spherical case at time = 1.0 ms. The numerical simulation is implemented with the code \textit{EBI-DNS}.
}\label{fig:tur_tri}
\end{figure*}

\paragraph{2D spherical case}
The computational domain of the spherical reaction-diffusion case, with size of $1{\rm cm} \times 1{\rm cm}$ and cells of 120,000, is filled with premixed gases characterized by a pressure of 1 atm and an equivalence ratio of 1. The ignition source is located at the center of the inlet with a radius of $0.4~\rm{mm}$. As shown in Fig.~\ref{fig:tur_tri}(A), DeePODE models can capture the flame structure and propagation for n-heptane. 
Moreover, Fig.\ref{fig:tur_tri}(C) shows that DeePODE can accurately predict the temperature profile distribution across the cross-section.

\paragraph{2D turbulent case} 

In the turbulent reaction-diffusion case, we set $512 \times 512$ cells for the computational domain of $1.5{\rm cm} \times 1.5 {\rm cm}$, the velocity field is generated by Passot-Pouquet isotropic kinetic energy spectrum. We set an initial field with $T_0= 300~{\rm K}$, $p_0 = 1~{\rm atm}$, $\phi_0 = 1$ in the domain. The ignition round is set in the center of the domain with a radius of $0.4{\rm mm}$. 
As illustrated in Fig.~\ref{fig:tur_tri}
(B) for the DRM19 mechanism, we compare the temperature distribution from 0.2 ms to 1.6 ms, obtained using CVODE and DeePODE. The DeePODE model demonstrates its capability to accurately capture the evolution of the turbulent reaction-diffusion system.

\paragraph{3D turbulent jet case} 
We evaluate our method using a well-established turbulent jet case as a benchmark for reaction-diffusion system simulation. This case provides comprehensive experimental and numerical data for thorough validation purposes.\cite{schneider2003flow,barlow1998effects}. 

The case features a main jet composed of ${\rm CH}_4$ and air, with an equivalence ratio of $\phi= 3.17$, a bulk velocity of $u_b= 49.6~{\rm m/s}$, and a temperature of $T= 294~{\rm K}$, delivered through a nozzle with a diameter of $7.2~{\rm mm}$. The flame is stabilized with a pilot jet of ${\rm C}_2{\rm H}_2$, ${\rm H}_2$, air, ${\rm CO}_2$, and ${\rm N}_2$, with a flow velocity of $11.4~{\rm m/s}$. The pilot nozzle has inner and outer diameters of $7.7~{\rm mm}$ and $18.2~{\rm mm}$, respectively. 
The jet inlets are surrounded by an air co-flow with an inner diameter of $18.9~{\rm mm}$. The computational domain contains 2,000,000 cells, and the chemical model employed is DRM19. 
Further setup details can be found in DLBFoam \cite{Morev2022, tekgul2021dlbfoam}. 

\begin{figure*}[htbp]
	\centering
        \subfigure{\includegraphics[width=1\textwidth]{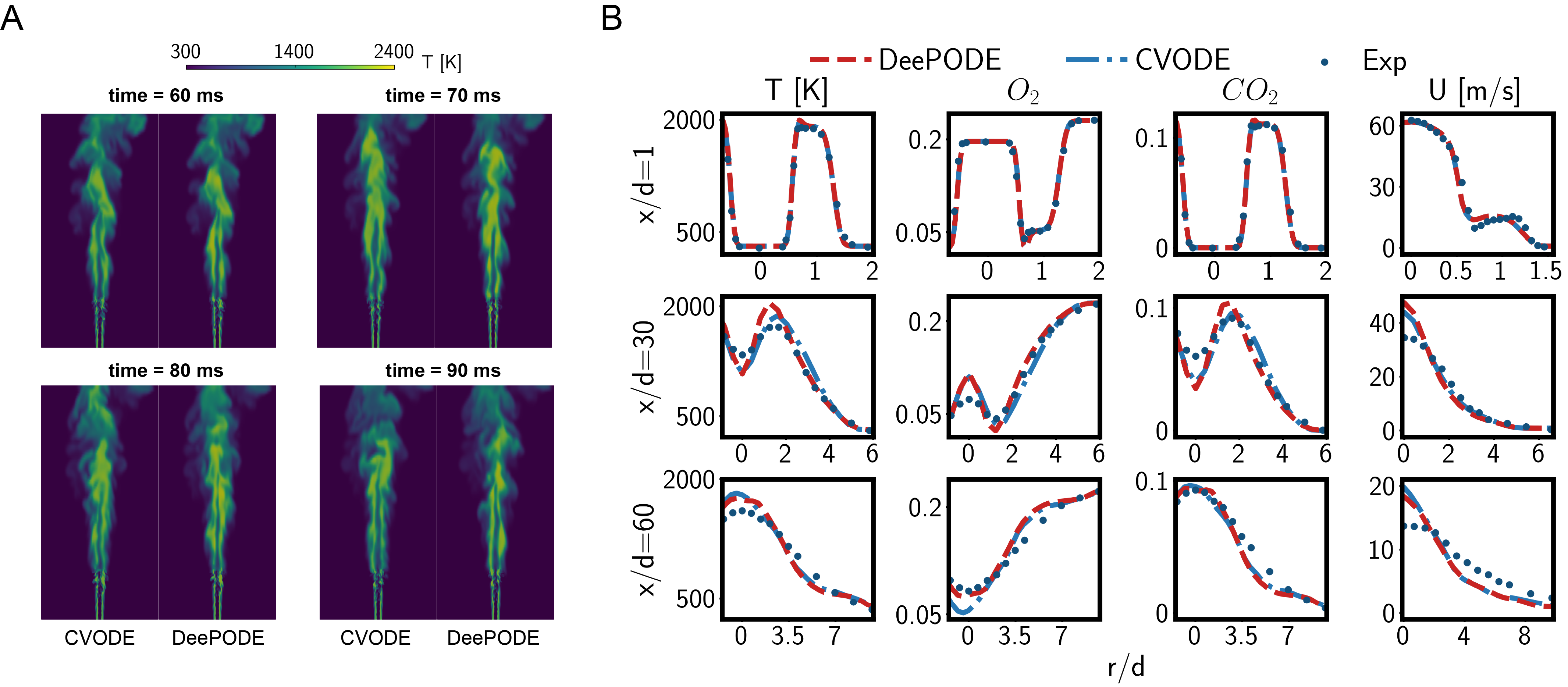}}
  \caption{\textbf{Three-dimensional turbulent jet case}. \textbf{(A)} shows the snapshots of temperature distribution by \textit{EBI-DNS} with CVODE or DeePODE for DRM19 at different time. \textbf{(B)} compares the average distributions of temperature, ${\rm O_2}$, ${\rm CO_2}$ and velocity by DeePODE, CVODE and experiment along radial direction on different axial locations .} \label{fig:sandia}
\end{figure*}

As shown in Fig. \ref{fig:sandia}(A), the simulation with DeePODE is stable and closely resembles the CVODE results across various time snapshots. 
The quantitative comparison among DeePODE, CVODE, and experimental data is presented in Fig. \ref{fig:sandia}(B), depicting the average radial distributions of temperature, ${\rm O}_2$, ${\rm CO}_2$, and velocity at different axial locations, specifically $x/d=1$, $30$, and $60$. The numerical and experimental result comparisons confirm DeePODE's exceptional performance in both accuracy and robustness. These findings reveal that data-driven methods can successfully handle traditionally challenging cases when provided with well-captured training data.

\subsection{Analysis}

\subsubsection{Error propagation}
\begin{figure*}[htbp]
    \centering
    \includegraphics[width=1\textwidth]{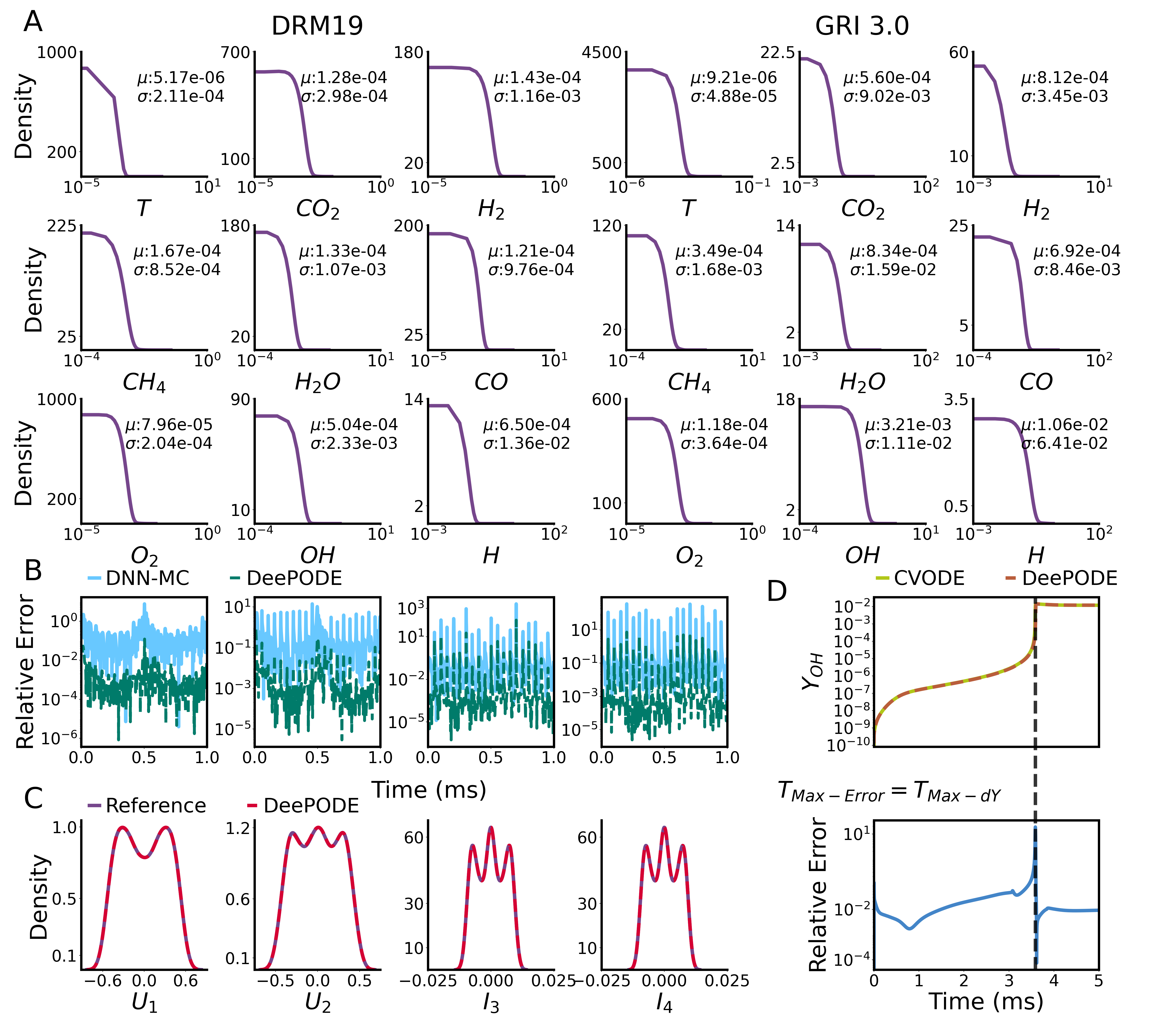}
    \caption{\textbf{Error analysis for DeePODE}.   \textbf{(A)} Distribution of relative error of one-step predictions in reaction-diffusion examples, on the 500,000 testing data from the zero- and one-dimensional manifold. \textbf{(B)} compares the error profiles of the EMCS and MC sampling methods in the context of an electronic dynamical simulation example. \textbf{(C)} compares the distribution of the evolution data obtained from the DeePODE approach with that of a traditional numerical solver. \textbf{(D)} illustrates the evolution trajectories predicted by DeePODE and CVODE, along with their absolute difference, in a reaction dynamics simulation example.}\label{fig:error_OH}
\end{figure*}

The error analysis presented in Fig. \ref{fig:error_OH} provides comprehensive insights into DeePODE's performance across different systems:
The empirical probability density function of relative error for one-step predictions in the reaction-diffusion system (Fig. \ref{fig:error_OH}(A)) demonstrates DeePODE's high accuracy across most test cases.
For electronic dynamics (Fig. \ref{fig:error_OH}(B)(C)), two key findings emerge. First, unlike traditional numerical schemes, DeePODE's error does not accumulate over the test period. Second, DeePODE's predictions closely match reference solutions, with discrepancies primarily manifesting as time-delay errors rather than fundamental prediction inaccuracies.
In reaction-diffusion simulations (Fig. \ref{fig:error_OH}(D)), the distribution of predicted quantities shows substantial overlap with benchmark solutions, confirming the method's reliability. While time-delay errors exist, they remain within acceptable bounds for practical applications, offering a favorable balance between computational efficiency and accuracy in complex system simulation.

\subsubsection{Computational efficiency}

We conduct a comprehensive performance analysis comparing DeePODE against numerical direct integration (DI), measuring the average computational time required per time step $\Delta t$. The evaluation is performed on two hardware platforms: a high-performance CPU (Intel® Xeon® Platinum 8260 CPU@2.40GHz) and a GPU accelerator (Tesla V100-SXM2-32GB). This dual-platform benchmarking provides insights into the computational efficiency of DeePODE across different computing architectures.

The computational performance data in Table \ref{tab:time_cost} demonstrates DeePODE's superior efficiency across multiple test cases. The method achieves speedup ratios of $10-100\times$ compared to traditional direct integration methods, showing impressive gains in both point-wise and parallel computing configurations. The most significant acceleration is observed in reaction-diffusion systems, where DeePODE exhibits its greatest performance advantages. These quantitative results validate DeePODE's ability to substantially improve computational efficiency while preserving solution quality, making it a powerful alternative to conventional approaches.

The superior computational performance of DeePODE can be attributed to two key advantages:
1) DeePODE exhibits a key advantage through its independence from system stiffness, which sets it apart from traditional integrators. ODE stiffness fundamentally limits traditional methods, forcing them to use extremely small time steps, particularly when handling high-dimensional systems. In contrast, DeePODE's stiffness-independent operation enables larger time steps while maintaining numerical stability. This independence from stiffness constraints represents a significant advancement in numerical integration methodology.
2) A key advantage of DeePODE lies in its inherent parallel processing capabilities. The method's architecture naturally performs matrix operations optimized for GPU acceleration. This feature is particularly beneficial for large-scale PDE-ODE coupled simulations, such as reaction-diffusion systems, where calculations must be performed across multiple grid points with local source terms. 

\renewcommand{\arraystretch}{1.0} 
\setlength{\tabcolsep}{0.2cm} 
\begin{table*}[htbp]
    \centering
    \begin{threeparttable}
        \begin{tabular}
        { @{}l @{\hspace{1pt}} c@{\hspace{2pt}}  c@{\hspace{4pt}} c@{\hspace{2pt}} c@{\hspace{2pt}} c@{\hspace{2pt}} c@{\hspace{2pt}} c@{\hspace{4pt}} c@{\hspace{2pt}} } 
            \toprule
            \multirow{2}{*}{\textbf{Test case}} & \multirow{2}{*}{\textbf{Model}} & \multirow{2}{*}{\textbf{\#dim}} & \multirow{2}{*}{\textbf{\#mesh}} & \multicolumn{3}{c}{\textbf{Time cost [s]~  $\downarrow$ }    } & \multicolumn{2}{c}{\textbf{Speed-up\tnote{1}~~$\uparrow$}}                                                               \\
                                                &                                 &                                 &                                  & \textbf{DI}                                                   & \textbf{DNN-CPU}                                & \textbf{DNN-GPU} & \textbf{CPU}               & \textbf{GPU}                  \\

            \midrule
            \rowcolor{lightblue} \multicolumn{2}{l}{ \textbf{Normal dynamical system}   }                                                                                                                                                                                                                                              \\
            Electronic dynamics                 & -                               & 15                              & 1                               & 0.03871                                                         & 0.00087                                                 & 0.00077                  &  \textbf{44.49}               &     \textbf{50.27}                   \\
            Electrolyte thermal runaway                & DMC                             & 104                            & 1                                & 0.0094                                                          &  0.001                                               &  0.0008                &  \textbf{9.4}                 &  \textbf{11.75}                   \\
            \midrule
            \rowcolor{lightpink} \multicolumn{2}{l}{ \textbf{Reaction-diffusion system} }                                                                                                                                                                                                                                              \\
            \text{1D laminar case}             & DRM19                           & 23                              & 4335                             & $1.75$                                                        & $0.178$                                         & $0.0034$         & \textbf{9.8}      & \textbf{51.29 }      \\
            \text{2D counterflow case}              & DRM19                           & 23                              & 960k                             & $174.57$                                                      & $32.53$                                         & $6.078$          & \textbf{5.37  }   & \textbf{28.72}       \\
            \text{2D spherical case}           & DRM19                           & 23                              & 120k                             & $50.10$                                                       & $4.17$                                          & $0.769$          & \textbf{12.01}    & \textbf{65.15}       \\
            \text{2D turbulent case}        & DRM19                           & 23                              & 262k                             & $48.04$                                                       & $11.36$                                         & $1.635$          & \textbf{4.23 }    & \textbf{29.38 }      \\
            \text{3D Sandia flame D\tnote{2} }            & DRM19                           & 23                              & 100k                             & $170.82$                                                      & $71.17$                                         & $1.687$          & \textbf{2.4 }     & \textbf{101.26}      \\
            \hline
            \text{1D laminar case}             & GRI 3.0                         & 55                              & 4335                             & $11.36$                                                       & $0.44$                                          & $0.042$          & \textbf{25.76 }   & \textbf{270.54}      \\
            \text{2D counterflow case}              & GRI 3.0                         & 55                              & 960k                             & $1228.97$                                                     & $31.98$                                         & $7.92$           & \textbf{38.43 }   & \textbf{155.23    }  \\
            \text{2D spherical case}           & GRI 3.0                         & 55                              & 120k                             & $270.02$                                                      & $4.93$                                          & $0.99$           & \textbf{54.75   } & \textbf{272.20   }   \\
            \hline
            \text{1D laminar case}             & n-heptane                       & 36                              & 4335                             & $3.92$                                                        & $0.25$                                          & $0.004$          & \textbf{16.29  }  & \textbf{98.05   }    \\
            \text{2D counterflow case}              & n-heptane                       & 36                              & 960k                             & $760.55$                                                      & $36.74$                                         & $6.502$          & \textbf{20.07 }   & \textbf{116.97 }     \\
            \text{2D spherical case}           & n-heptane                       & 36                              & 120k                             & $204.07$                                                      & $4.89$                                          & $0.825$          & \textbf{ 41.70}   & \textbf{247.36}      \\
            \bottomrule
        \end{tabular}
        \begin{tablenotes}
			\footnotesize
			\item[1] Compared with numerical direct integrator (DI) on CPU. 
                \item[2] A turbulent reaction-diffusion benchmark case. Only consider the 100k grids that DNN applies for.
		\end{tablenotes}
		\caption{Time cost comparison.}
        \label{tab:time_cost}
    \end{threeparttable}
\end{table*}

\subsubsection{Characteristic time analysis}
\label{sec:csp_time_gms}

The temporal evolution of randomly generated data plays a vital role in capturing representative dynamics in high-dimensional phase space. We conduct a comparative analysis of different evolution strategies and their associated characteristic timescales to demonstrate this importance.

Using the 22-dimensional DRM19 chemical reaction model as an example, we investigates three distinct strategies for selecting $\boldsymbol{\tau}$: fixed $\boldsymbol{\tau}$ using constant evolution time intervals, fully adaptive $\boldsymbol{\tau}$ dynamically adjusting evolution times based on system behavior, and increasing sequence $\boldsymbol{\tau}$ with progressively longer evolution times. This comparison helps us understand how different time evolution strategies affect the method's performance.
The fixed $\boldsymbol{\tau}$ strategy uses constant time intervals of $\tau_i = 10^{-6}$ s, with evolution steps $k$ set to either 10 or 20. The adaptive $\boldsymbol{\tau}$ approach, tested with $k=10$ and $k=20$, determines each $\tau_i$ based on the local state vector's integration timescale. We employ the computational singular perturbation (CSP) \cite{lam1994csp} method to identify the intrinsic timescale $\tau_{csp}$, derived from the derivative of the Jacobian matrix. This choice is justified by $\tau_{csp}$'s ability to represent local timescales, with detailed CSP parameters in Appendix \ref{sec:csp}. The increasing sequence strategy uses $k=10$ with progressively longer evolution times: $\boldsymbol{\tau}=[10^{-6}, 10^{-6}, 10^{-5}, 10^{-5}, 10^{-4}, 10^{-4}, 10^{-3}, 2\times10^{-3}, 5\times 10^{-3}]$ s. All strategies are initialized with the same Monte Carlo dataset to ensure a fair comparison.

\begin{figure}[htbp]
	\centering
        \includegraphics[width=1\textwidth]{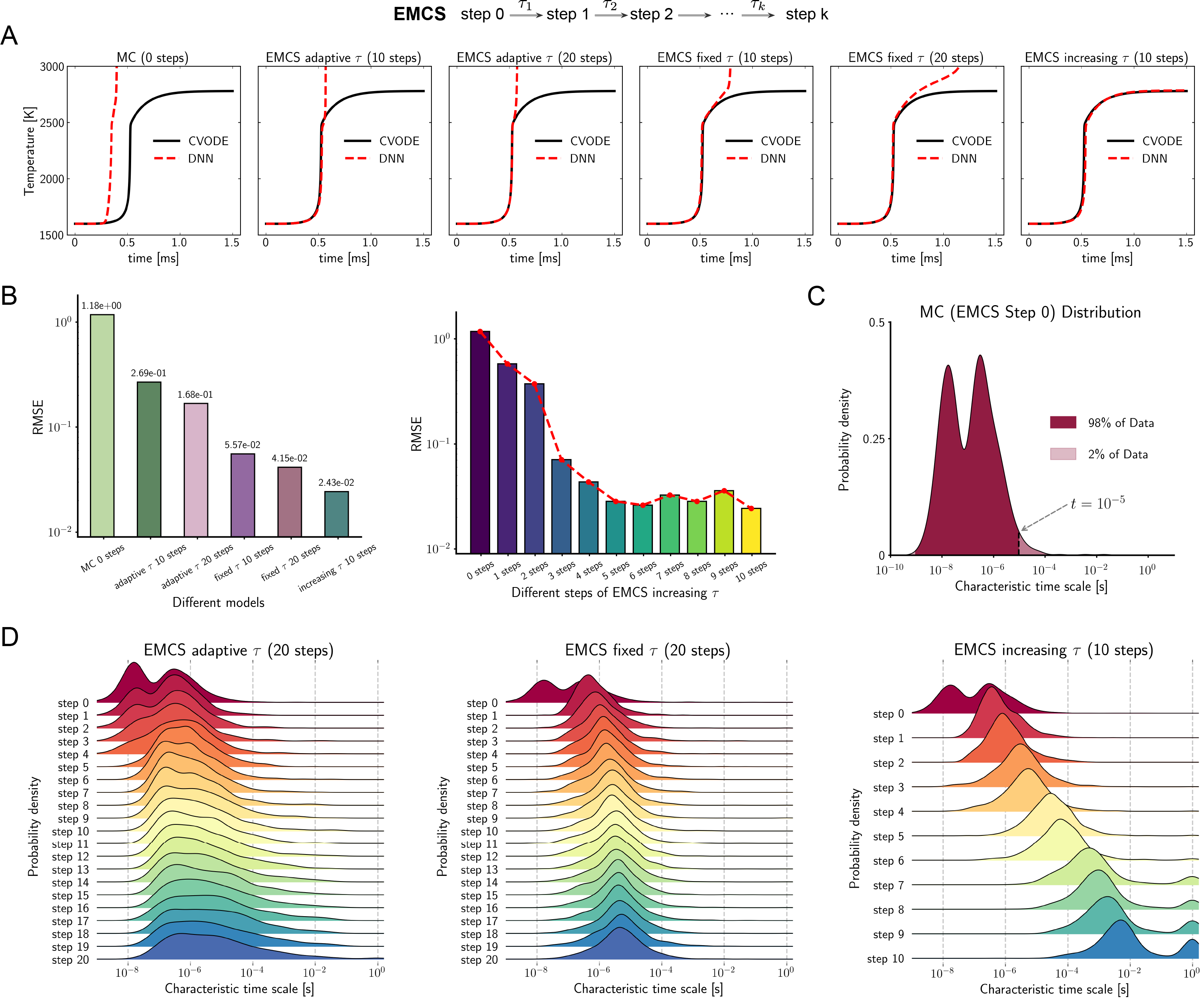}
\caption{\textbf{
Different evolution time selection strategies for EMCS}. \textbf{(A)} shows temperature evolution predictions from six DNN models compared against ground truth solutions (black solid lines) from the CVODE solver. Initial conditions are set at 1 atm pressure, 1600 K temperature, and stoichiometric equivalence ratio, with a total evolution time of 1.5 ms. \textbf{(B)} quantifies prediction accuracy through temperature RMSE. The left plot compares six models: basic Monte Carlo (no evolution), adaptive $\boldsymbol{\tau}$ (10 and 20 steps), fixed $\boldsymbol{\tau}$ (10 and 20 steps), and increasing $\boldsymbol{\tau}$ (10 steps). The right plot examines how RMSE varies with evolution steps ($k=0$ to 10) for the increasing $\boldsymbol{\tau}$ implementation. \textbf{(C)} illustrates the characteristic timescale distribution of the baseline Monte Carlo dataset (step 0 of EMCS). \textbf{(D)} compares how characteristic timescale distributions evolve across different steps for all three EMCS strategies.
} 
\label{fig:csp_time_steps}
\end{figure}

Our comparative analysis examines six models trained on different EMCS datasets. As shown in Fig. \ref{fig:csp_time_steps}(A), the traditional Monte Carlo-based DNN model poorly predicts temperature trajectories. In contrast, EMCS-based models using fixed and adaptive time sequences (with 10 and 20 steps) demonstrated moderate improvements. Most notably, the EMCS model employing increasing time intervals achieves superior prediction accuracy. These qualitative observations in Fig. \ref{fig:csp_time_steps}(B) are supported by quantitative analysis of point-to-point temperature prediction RMSE, where the increasing-interval EMCS model achieved the lowest error rate of 0.0243. Additionally, we found that prediction accuracy improves systematically with the number of evolution steps, as evidenced by decreasing RMSE values across evolution steps for $k$ ranging from 0 to 10.

Analysis of characteristic timescales reveals why EMCS excels at handling multi-scale dynamics. Using $\tau_{csp}$ to measure characteristic time in the data, we find that traditional Monte Carlo sampling exhibits a significant limitation: it concentrates 98\% of data points below $10^{-5}$ seconds, failing to adequately represent the system's multiple timescales as shown in Fig. \ref{fig:csp_time_steps}(C).
In contrast, Fig. \ref{fig:csp_time_steps}(D) demonstrates how EMCS improves through its evolution steps, with each step covering a distinct range of timescales. The findings reveal a correlation between DNN performance and the breadth of timescale coverage in the training data. The three EMCS strategies demonstrate in the figure show that broader coverage of local timescales leads to improved prediction accuracy. 
The observation reveals a fundamental insight: comprehensive timescale coverage is critical in enhancing DNN's ability to handle transient multiscale dynamical systems. This explains why proper sampling across different timescales directly translates to improved prediction accuracy, suggesting that the quality of multiscale predictions is fundamentally linked to how well the training data spans the relevant temporal scales. This systematic evolution of characteristic time distribution across steps is a consistent pattern observed across various examples, as illustrated in Fig. \ref{fig:csp_time_preadtor_modulator}.

\section{Discussion}

We introduce DeePODE, a novel deep-learning approach for solving high-dimensional multiscale dynamical systems. The method combines two key innovations: an Evolutionary Monte Carlo Sampling (EMCS) technique and a specialized deep neural network architecture. EMCS generates representative training samples that effectively capture the multiple timescales inherent in complex dynamical systems. The neural network design incorporates large time steps to overcome system stiffness limitations, substantially improving computational efficiency. This integrated framework enables accurate simulation of complex dynamical systems while achieving significant speedup compared to traditional numerical methods.

The DeePODE method is comprehensively validated in diverse multiscale systems across various fields, such as the predator-prey model in the ecosystem, the ring modulator model in the power system, an electrolyte thermal runaway in the lithium battery design, and a series of reaction-diffusion simulations considering detailed chemical kinetics. 
In particular, for the 15-dimensional ring modulator and 104-dimensional electrolyte thermal runaway cases, 
we compare the performance of DNN-MC/DeePODE with traditional implicit ODE solvers. 
The results reveal that the DeePODE method achieves the accuracy of implicit numerical schemes while maintaining the computational speed of explicit schemes, offering a robust and efficient alternative for solving multiscale dynamical systems. 

Furthermore, we evaluate DeePODE's performance in complex flow scenarios governed by partial differential equations through extensive benchmark testing. The method accurately predicts reaction-diffusion system dynamics while achieving computational acceleration of one to two orders of magnitude over conventional numerical solvers. A notable strength of DeePODE lies in its generalization capability. Despite using training data not directly derived from specific simulations, the method accurately predicts turbulent reaction front behavior under unseen conditions. This robust generalization can be attributed to the comprehensive multiscale sampling achieved through EMCS, which captures essential system dynamics across different temporal and spatial scales.

The superior performance of EMCS in addressing multiscale dynamical systems arises from its innovative approach to adaptive data sampling and temporal scale coverage. Through its evolutionary process, EMCS adapts to the local gradient structure of the phase space, efficiently allocating samples where they are most needed. This approach generates fewer data points in regions of low gradient—where system behavior is more predictable—aligning with neural networks' natural ability to learn smooth, low-frequency patterns. As a result, DeePODE achieves accurate approximations across diverse operating conditions while maintaining an economical dataset size. Another distinctive strength of EMCS is its comprehensive coverage of characteristic timescales in high-dimensional systems, leading to robust generalization in temporal evolution prediction. As evidenced in Fig.~\ref{fig:csp_time_steps}(D), EMCS samples effectively across the full spectrum of system dynamics, whereas traditional Monte Carlo methods fail to capture this temporal diversity. This thorough representation of multiple timescales fundamentally enhances the method's predictive accuracy and reliability. \red{It's noted that in deep learning methodologies, dataset size typically has a significant influence on model performance. Under appropriate conditions, larger training datasets generally yield improved model outcomes – a principle that equally applies to the EMCS sampling approach.}

In summary, DeePODE demonstrates exceptional capability in simulating complex dynamical systems, with promising applications ranging from atmospheric pollution forecasting~\cite{chae2021pm10} to biological oscillator modulation~\cite{Linwei2023}. Its combination of efficiency and accuracy makes it a valuable tool for addressing challenging multiscale problems in science and engineering.

\section*{Declaration of interests}
The authors declare that they have no known competing financial interests or personal relationships that could have appeared to influence the work reported in this paper.

\section*{Declaration of generative AI in scientific writing}

During the preparation of this work the authors used generative AI in order to polish the language. After using this tool, the authors reviewed and edited the content as needed and take full responsibility for the content of the publication.

\section*{Acknowledgements}

This work is sponsored by the National Natural Science Foundation of China 92470127 (T. Z.), No. 62002221 (Z. X.), 12371511 (Z. X.), 12422119 (Z. X.), 92270203 (T. Z.), Grant No. 12101402 (Y. Z.), the Lingang Laboratory Grant No.LG-QS-202202-08 (Y. Z.), Shanghai Municipal of Science and Technology Project Grant No. 20JC1419500 (Y. Z.), Shanghai Municipal of Science and Technology Major Project No. 2021SHZDZX0102, and the HPC of School of Mathematical Sciences and the Student Innovation Center, and the Siyuan-1 cluster supported by the Center for High-Performance Computing at Shanghai Jiao Tong University, Key Laboratory of Marine Intelligent Equipment and System, Ministry of Education, P.R. China. SJTU Kunpeng \& Ascend Center of Excellence partially supported this work.

\bibliographystyle{elsarticle-num}
\bibliography{references} 

\begin{thebibliography}{10}
\expandafter\ifx\csname url\endcsname\relax
  \def\url#1{\texttt{#1}}\fi
\expandafter\ifx\csname urlprefix\endcsname\relax\def\urlprefix{URL }\fi
\expandafter\ifx\csname href\endcsname\relax
  \def\href#1#2{#2} \def\path#1{#1}\fi

\bibitem{barlow1998effects}
R.~Barlow, J.~Frank, Effects of turbulence on species mass fractions in methane/air jet flames, in: Symposium (International) on Combustion, Vol.~27, Elsevier, 1998, pp. 1087--1095.

\bibitem{mcadams1995circuit}
H.~H. McAdams, L.~Shapiro, Circuit simulation of genetic networks, Science 269~(5224) (1995) 650--656.

\bibitem{volterra1931variations}
V.~Volterra, Variations and fluctuations of the number of individuals in animal species living together., Animal ecology (1931) 412--433.

\bibitem{ye2016analysis}
H.~Ye, Y.~Liu, P.~Zhang, Z.~Du, Analysis and detection of forced oscillation in power system, IEEE Transactions on Power Systems 32~(2) (2016) 1149--1160.

\bibitem{sun2009library}
Y.~Sun, D.~Zhou, A.~V. Rangan, D.~Cai, Library-based numerical reduction of the hodgkin--huxley neuron for network simulation, Journal of computational neuroscience 27 (2009) 369--390.

\bibitem{dissanayake1994neural}
M.~Dissanayake, N.~Phan-Thien, Neural-network-based approximations for solving partial differential equations, communications in Numerical Methods in Engineering 10~(3) (1994) 195--201.

\bibitem{raissi2019physics}
M.~Raissi, P.~Perdikaris, G.~E. Karniadakis, Physics-informed neural networks: A deep learning framework for solving forward and inverse problems involving nonlinear partial differential equations, Journal of Computational Physics 378 (2019) 686--707.

\bibitem{zhang2022multi}
T.~Zhang, Y.~Yi, Y.~Xu, Z.~X. Chen, Y.~Zhang, W.~E, Z.-Q.~J. Xu, A multi-scale sampling method for accurate and robust deep neural network to predict combustion chemical kinetics, Combustion and Flame 245 (2022) 112319.
\newblock \href {https://doi.org/https://doi.org/10.1016/j.combustflame.2022.112319} {\path{doi:https://doi.org/10.1016/j.combustflame.2022.112319}}.

\bibitem{li2023mathematical}
S.~Li, D.~W. McLaughlin, D.~Zhou, Mathematical modeling and analysis of spatial neuron dynamics: Dendritic integration and beyond, Communications on Pure and Applied Mathematics 76~(1) (2023) 114--162.

\bibitem{LED2022}
P.~R. Vlachas, G.~Arampatzis, C.~Uhler, P.~Koumoutsakos, \href{https://doi.org/10.1038/s42256-022-00464-w}{Multiscale simulations of complex systems by learning their effective dynamics}, Nature Machine Intelligence 4~(4) (2022) 359--366.
\newblock \href {https://doi.org/10.1038/s42256-022-00464-w} {\path{doi:10.1038/s42256-022-00464-w}}.
\newline\urlprefix\url{https://doi.org/10.1038/s42256-022-00464-w}

\bibitem{AdaLED2023}
I.~Kičić, P.~R. Vlachas, G.~Arampatzis, M.~Chatzimanolakis, L.~Guibas, P.~Koumoutsakos, \href{http://dx.doi.org/10.1016/j.cma.2023.116204}{2023 adaled adaptive learning of effective dynamics for online modeling of complex systems}, Computer Methods in Applied Mechanics and Engineering (2023) 116204\href {https://doi.org/10.1016/j.cma.2023.116204} {\path{doi:10.1016/j.cma.2023.116204}}.
\newline\urlprefix\url{http://dx.doi.org/10.1016/j.cma.2023.116204}

\bibitem{GLED2024}
H.~Gao, S.~Kaltenbach, P.~Koumoutsakos, Generative learning for forecasting the dynamics of high-dimensional complex systems, Nature Communications 15~(1) (Oct 2024).
\newblock \href {https://doi.org/10.1038/s41467-024-53165-w} {\path{doi:10.1038/s41467-024-53165-w}}.

\bibitem{milman2009asymptotic}
V.~D. Milman, G.~Schechtman, Asymptotic theory of finite dimensional normed spaces: Isoperimetric inequalities in riemannian manifolds, Vol. 1200, 2009.

\bibitem{chen1989pdf}
J.-Y. Chen, W.~Kollmann, R.~Dibble, Pdf modeling of turbulent nonpremixed methane jet flames, Combustion Science and Technology 64~(4-6) (1989) 315--346.

\bibitem{pope1997computationally}
S.~B. Pope, Computationally efficient implementation of combustion chemistry using in situ adaptive tabulation (1997).

\bibitem{christo1996artificial}
F.~Christo, A.~Masri, E.~Nebot, Artificial neural network implementation of chemistry with pdf simulation of h2/co2 flames, Combustion and Flame 106~(4) (1996) 406--427.

\bibitem{blasco1998modelling}
J.~Blasco, N.~Fueyo, C.~Dopazo, J.~Ballester, Modelling the temporal evolution of a reduced combustion chemical system with an artificial neural network, Combustion and Flame 113~(1-2) (1998) 38--52.

\bibitem{choi2005fast}
Y.~Choi, J.-Y. Chen, Fast prediction of start-of-combustion in hcci with combined artificial neural networks and ignition delay model, Proceedings of the Combustion Institute 30~(2) (2005) 2711--2718.

\bibitem{sen2010linear}
B.~A. Sen, S.~Menon, Linear eddy mixing based tabulation and artificial neural networks for large eddy simulations of turbulent flames, Combustion and Flame 157~(1) (2010) 62--74.

\bibitem{sinaei2017large}
P.~Sinaei, S.~Tabejamaat, Large eddy simulation of methane diffusion jet flame with representation of chemical kinetics using artificial neural network, Proceedings of the Institution of Mechanical Engineers, Part E: Journal of Process Mechanical Engineering 231~(2) (2017) 147--163.

\bibitem{zhang2021deep}
T.~Zhang, Y.~Zhang, W.~E, Y.~Ju, A deep learning-based ode solver for chemical kinetics, in: AIAA Science and Technology Forum and Exposition, AIAA SciTech Forum 2021, American Institute of Aeronautics and Astronautics Inc, AIAA, 2021.

\bibitem{owoyele2022chemnode}
O.~Owoyele, P.~Pal, Chemnode: A neural ordinary differential equations framework for efficient chemical kinetic solvers, Energy and AI 7 (2022) 100118.

\bibitem{almeldein2022accelerating}
A.~Almeldein, N.~Van~Dam, Accelerating chemical kinetics calculations with physics informed neural networks, in: Internal Combustion Engine Division Fall Technical Conference, Vol. 86540, American Society of Mechanical Engineers, 2022, p. V001T06A007.

\bibitem{de2022physics}
M.~De~Florio, E.~Schiassi, R.~Furfaro, Physics-informed neural networks and functional interpolation for stiff chemical kinetics, Chaos: An Interdisciplinary Journal of Nonlinear Science 32~(6) (2022) 063107.

\bibitem{yao2022gradient}
S.~Yao, A.~Kronenburg, A.~Shamooni, O.~Stein, W.~Zhang, Gradient boosted decision trees for combustion chemistry integration, Applications in Energy and Combustion Science 11 (2022) 100077.

\bibitem{silver2016mastering}
D.~Silver, A.~Huang, C.~J. Maddison, A.~Guez, L.~Sifre, G.~Van Den~Driessche, J.~Schrittwieser, I.~Antonoglou, V.~Panneershelvam, M.~Lanctot, et~al., Mastering the game of go with deep neural networks and tree search, nature 529~(7587) (2016) 484--489.

\bibitem{hawkes2005direct}
E.~R. Hawkes, R.~Sankaran, J.~C. Sutherland, J.~H. Chen, Direct numerical simulation of turbulent combustion: fundamental insights towards predictive models, in: Journal of Physics: Conference Series, Vol.~16, IOP Publishing, 2005, p.~65.

\bibitem{mettler2012top}
M.~S. Mettler, D.~G. Vlachos, P.~J. Dauenhauer, Top ten fundamental challenges of biomass pyrolysis for biofuels, Energy \& Environmental Science 5~(7) (2012) 7797--7809.

\bibitem{seinfeld1998air}
J.~H. Seinfeld, S.~N. Pandis, From air pollution to climate change, Atmospheric chemistry and physics 1326 (1998).

\bibitem{Box1964}
G.~E. Box, D.~R. Cox, An analysis of transformations, Journal of the Royal Statistical Society: Series B (Methodological) 26~(2) (1964) 211--243.

\bibitem{xu_training_2018}
Z.-Q.~J. Xu, Y.~Zhang, Y.~Xiao, Training behavior of deep neural network in frequency domain, International Conference on Neural Information Processing (2019) 264--274.

\bibitem{xu2019frequency}
Z.-Q.~J. Xu, Y.~Zhang, T.~Luo, Y.~Xiao, Z.~Ma, Frequency principle: Fourier analysis sheds light on deep neural networks, Communications in Computational Physics 28~(5) (2020) 1746--1767.

\bibitem{qin2021data}
T.~Qin, Z.~Chen, J.~D. Jakeman, D.~Xiu, Data-driven learning of nonautonomous systems, SIAM Journal on Scientific Computing 43~(3) (2021) A1607--A1624.

\bibitem{DeePMR2022}
Z.~Wang, Y.~Zhang, E.~Zhao, Y.~Ju, W.~E, Z.-Q.~J. Xu, T.~Zhang, A deep learning-based model reduction (deepmr) method for simplifying chemical kinetics (2022).
\newblock \href {http://arxiv.org/abs/2201.02025} {\path{arXiv:2201.02025}}.

\bibitem{zirwes2020quasi}
T.~Zirwes, F.~Zhang, P.~Habisreuther, M.~Hansinger, H.~Bockhorn, M.~Pfitzner, D.~Trimis, Quasi-dns dataset of a piloted flame with inhomogeneous inlet conditions, Flow, Turbulence and Combustion 104~(4) (2020) 997--1027.

\bibitem{chen2009effects}
Z.~Chen, M.~P. Burke, Y.~Ju, Effects of lewis number and ignition energy on the determination of laminar flame speed using propagating spherical flames, Proceedings of the Combustion Institute 32~(1) (2009) 1253--1260.

\bibitem{zhang2021studies}
T.~Zhang, A.~J. Susa, R.~K. Hanson, Y.~Ju, Studies of the dynamics of autoignition assisted outwardly propagating spherical cool and double flames under shock-tube conditions, Proceedings of the Combustion Institute 38~(2) (2021) 2275--2283.

\bibitem{mao2022deepflame}
R.~Mao, M.~Lin, Y.~Zhang, T.~Zhang, Z.-Q.~J. Xu, Z.~X. Chen, Deepflame: A deep learning empowered open-source platform for reacting flow simulations, arXiv preprint arXiv:2210.07094 (2022).

\bibitem{cohen1996cvode}
S.~D. Cohen, A.~C. Hindmarsh, P.~F. Dubois, et~al., Cvode, a stiff/nonstiff ode solver in c, Computers in physics 10~(2) (1996) 138--143.

\bibitem{schneider2003flow}
C.~Schneider, A.~Dreizler, J.~Janicka, E.~Hassel, Flow field measurements of stable and locally extinguishing hydrocarbon-fuelled jet flames, Combustion and Flame 135~(1-2) (2003) 185--190.

\bibitem{Morev2022}
I.~Morev, B.~Tekg\"{u}l, M.~Gadalla, A.~Shahanaghi, J.~Kannan, S.~Karimkashi, O.~Kaario, V.~Vuorinen, \href{https://doi.org/10.1063/5.0077437}{Fast reactive flow simulations using analytical jacobian and dynamic load balancing in {OpenFOAM}}, Physics of Fluids 34~(2) (2022) 021801.
\newblock \href {https://doi.org/10.1063/5.0077437} {\path{doi:10.1063/5.0077437}}.
\newline\urlprefix\url{https://doi.org/10.1063/5.0077437}

\bibitem{tekgul2021dlbfoam}
B.~Tekg{\"u}l, P.~Peltonen, H.~Kahila, O.~Kaario, V.~Vuorinen, Dlbfoam: An open-source dynamic load balancing model for fast reacting flow simulations in openfoam, Computer Physics Communications (2021) 108073.

\bibitem{lam1994csp}
S.~Lam, D.~Goussis, The csp method for simplifying kinetics, International journal of chemical kinetics 26~(4) (1994) 461--486.

\bibitem{chae2021pm10}
S.~Chae, J.~Shin, S.~Kwon, S.~Lee, S.~Kang, D.~Lee, Pm10 and pm2.5 real-time prediction models using an interpolated convolutional neural network, Scientific Reports 11~(1) (2021) 11952.
\newblock \href {https://doi.org/10.1038/s41598-021-93696-4} {\path{doi:10.1038/s41598-021-93696-4}}.

\bibitem{Linwei2023}
K.~Wang, L.~Yang, S.~Zhou, W.~Lin, \href{https://doi.org/10.1063/5.0167555}{Desynchronizing oscillators coupled in multi-cluster networks through adaptively controlling partial networks}, Chaos: An Interdisciplinary Journal of Nonlinear Science 33~(9) (2023) 091101.
\newblock \href {https://doi.org/10.1063/5.0167555} {\path{doi:10.1063/5.0167555}}.
\newline\urlprefix\url{https://doi.org/10.1063/5.0167555}

\bibitem{maas1992simplifying}
U.~Maas, S.~B. Pope, Simplifying chemical kinetics: intrinsic low-dimensional manifolds in composition space, Combustion and flame 88~(3-4) (1992) 239--264.

\bibitem{Timofeev2022modulator}
K.~A. Timofeev, Y.~V. Shornikov, Computer simulation of dynamic processes in electronic systems, in: 2022 IEEE International Multi-Conference on Engineering, Computer and Information Sciences (SIBIRCON), 2022, pp. 1860--1863.
\newblock \href {https://doi.org/10.1109/SIBIRCON56155.2022.10016988} {\path{doi:10.1109/SIBIRCON56155.2022.10016988}}.

\bibitem{chi2021fly}
C.~Chi, G.~Janiga, D.~Th{\'e}venin, On-the-fly artificial neural network for chemical kinetics in direct numerical simulations of premixed combustion, Combustion and Flame 226 (2021) 467--477.

\bibitem{chaos2007high}
M.~Chaos, A.~Kazakov, Z.~Zhao, F.~L. Dryer, A high-temperature chemical kinetic model for primary reference fuels, International Journal of Chemical Kinetics 39~(7) (2007) 399--414.

\bibitem{weng2025CF}
Y.~Weng, H.~Li, H.~Zhang, Z.~X. Chen, D.~Zhou, Extended fourier neural operators to learn stiff chemical kinetics under unseen conditions, Combustion and Flame 272 (2025) 113847.
\newblock \href {https://doi.org/10.1016/j.combustflame.2024.113847} {\path{doi:10.1016/j.combustflame.2024.113847}}.

\bibitem{goswami2024CMAM}
S.~Goswami, A.~D. Jagtap, H.~Babaee, B.~T. Susi, G.~E. Karniadakis, Learning stiff chemical kinetics using extended deep neural operators, Computer Methods in Applied Mechanics and Engineering 419 (2024) 116674.
\newblock \href {https://doi.org/10.1016/j.cma.2023.116674} {\path{doi:10.1016/j.cma.2023.116674}}.

\bibitem{ding2021CF}
T.~Ding, T.~Readshaw, S.~Rigopoulos, W.~Jones, Machine learning tabulation of thermochemistry in turbulent combustion: An approach based on hybrid flamelet/random data and multiple multilayer perceptrons, Combustion and Flame 231 (2021) 111493.
\newblock \href {https://doi.org/10.1016/j.combustflame.2021.111493} {\path{doi:10.1016/j.combustflame.2021.111493}}.

\bibitem{readshaw2023CF}
T.~Readshaw, L.~L.~C. Franke, W.~P. Jones, S.~Rigopoulos, Simulation of turbulent premixed flames with machine learning - tabulated thermochemistry, Combustion and Flame 258 (2023) 113058.
\newblock \href {https://doi.org/10.1016/j.combustflame.2023.113058} {\path{doi:10.1016/j.combustflame.2023.113058}}.

\end{thebibliography}

\newpage

\begin{appendices}

\part{}
\vspace{-10em}

\setcounter{parttocdepth}{3}
\part{Appendix} 
\parttoc 

\section{Additional Method Description} 

\subsection{EMCS Method}
Given an ODE system with the form
\begin{align}
        &\frac{d\boldsymbol{x}}{dt} = \boldsymbol{f}(\boldsymbol{x},t)\\
        &\boldsymbol{x}|_{t=0}=\boldsymbol{x}_0
\end{align}
where $\boldsymbol{x}(t) = [x_1(t),x_2(t),\dots,x_n(t)]$. EMCS for such a dynamics system can be formulated as

\begin{algorithm}[!ht]
            \renewcommand{\algorithmicrequire}{\textbf{Input:}}
    	\renewcommand{\algorithmicensure}{\textbf{Output:}}
    \caption{EMCS method}
    \label{alg:emcs}
    \begin{algorithmic}[1]
    \REQUIRE  System dimension $n$, the data size $N_{0}$ of MC sampling, the evolution time selection strategy $\boldsymbol{\tau} = [\tau_1, \tau_2, \cdots, \tau_k]$, the time step for DNN prediction $\Delta t$. 
    
    \STATE    \blue{ // Step 1: Range Estimation }
    \STATE Collect data points from multiple temporal trajectories and form the set $\boldsymbol{x}_{\rm MF}$;
    \STATE Compute the label set $\boldsymbol{u}_{\rm MF}$;
    \STATE Update the range estimation $R_{\rm MF}=[\min \boldsymbol{x}_{\rm MF},\max \boldsymbol{x}_{\rm MF}]\times [\min \boldsymbol{u}_{\rm MF},\max \boldsymbol{u}_{\rm MF}]$;

    \STATE    \blue{ // Step 2: Monte Carlo (MC) Sampling }
    \FOR{$i = 1 : n$}
    \IF{ the component $x^i$ spans multiple scales } 
        \STATE Sample $N_0$ sub-dimension samples $\log x^i\sim \mathcal{U}(\log\min \boldsymbol{x}^i_{\rm MF}, \log \max \boldsymbol{x}^i_{\rm MF})$;
    \ELSE 
        \STATE Sample $N_0$ sub-dimension samples $x^i\sim \mathcal{U}(\min \boldsymbol{x}^i_{\rm MF}, \max \boldsymbol{x}^i_{\rm MF})$;
    \ENDIF
    \ENDFOR
    \STATE Merge all the $N_0$ data samples into the set  $ \mathbb{D}_{\text{MC}}$.

    \STATE    \blue{ // Step 3: Evolution Augmented Generation}
    \FOR{ each $\boldsymbol{x}(t) \in \mathbb{D}_{\text{MC}} $}
        \STATE $\boldsymbol{x}_1:= \boldsymbol{x}(t)$;
        \FOR{$i = 1 : k$}
            \STATE $\boldsymbol{x}_{i+1}: = \boldsymbol{x}_i + \int_{t}^{t+\tau_i} \boldsymbol{f}(\boldsymbol{x}, s)\mathrm{d} s$.  \blue{ $\quad\triangleright~ \text{temporal evolution} $};
        \ENDFOR
        \FOR{$i = 1 : k+1$}
            \STATE $\boldsymbol{u}_{i} :=  \int_{t}^{t+\Delta t} \boldsymbol{f}(\boldsymbol{x}, s)\mathrm{d} s$;    \blue{ $\quad\triangleright~ \text{label generation} $ }
            \STATE Update the dataset $\mathbb{D}_{\text{EMCS}} = \mathbb{D}_{\text{EMCS}}\cup \{(\boldsymbol{x}_i, \boldsymbol{u}_i)\} \cap R_{\rm MF}  $.
    
         \ENDFOR
    \ENDFOR
    
    \ENSURE The final dataset $\mathbb{D}_{\text{EMCS}}$.
    \end{algorithmic}
    \end{algorithm}

\subsection{EMCS Method for Chemical Kinetics}\label{app:chemical}

Given a specific reactor and different initial temperatures, pressures, species concentrations, the trajectories of all species form a low-dimensional manifold \cite{maas1992simplifying}. Manifold sampling comprises two parts. First, we sample from zero-dimensional homogeneous reactions, where heated, homogeneous fuel and oxygen ignite spontaneously. The initial temperature ($T$) range is set so that the ignition delay time belongs to $[0.1 {\rm ms},10 {\rm ms}]$. Initial pressure $P$ belongs to $[0.5 {\rm atm},2 {\rm atm}]$. Equivalence ratio $\phi$ belongs to $[0.5,3]$. We use the open-source chemical kinetics library Cantera\footnote{\url{https://cantera.org/}} to run approximately 5000 cases with random initial conditions and maximal time step size $10^{-8}~{\rm s}$ until maximal running time (e.g., $10 {\rm ms}$) or when the temperature change is smaller than 0.001 within $10^{-7}~{\rm s}$.  We sample each state $\boldsymbol{x}(t)=[T(t),P(t),Y_1(t),\cdots,Y_{n_s}(t)]$ ($Y_i\in[0,1]$ is the mass fraction of species $i$) on each trajectory every $10^{\scriptscriptstyle{-7}}{\rm s}$ and its corresponding label $\boldsymbol{u}(t)=\boldsymbol{x}(t+\Delta t)-\boldsymbol{x}(t)$, where $\Delta t=10^{-6}  {\rm s}$ is the step size we use for DNN model. The other is sampling from one-dimensional freely propagating premixed flame example at 1 atm, $\phi=1$, and temperature is set from 300K to 1000K (every 10K). The physical length is 0.03 m  with 800 grid points. We sample from the stable states of grids. Both zero-dimensional and one-dimensional data constitute the manifold sample set, which is denoted as $\mathbb{D}_{\rm Init}=\{(\boldsymbol{x}(t), \boldsymbol{u}(t))\}$.

Next, we use the manifold data $D_{\rm Init}$ to filter data generated through global multiscale manifold sampling as follows. 
Next, we perform a global Monte-Carlo sampling. The pressure range is the same as for manifold sampling. The ranges of temperature and mass fractions of species are determined by the corresponding range in $\mathbb{D}_{\rm MF}$. Each data point $\boldsymbol{x}(t)$ is sampled uniformly in the log scale within the given range by the Monte-Carlo method with a sample size of approximately 3,200,000. The label for $\boldsymbol{x}(t)$ is computed by Cantera as $\boldsymbol{u}(t)=\boldsymbol{x}(t+\Delta t)-\boldsymbol{x}(t)$. Each data point is filtered by the manifold data as follows. A data point is kept if each element of its label belongs to a range determined by the manifold data, i.e., $u(t)_i\in[\lambda_1u_{{\rm min},i},\lambda_2u_{{\rm max},i}]$, where $u_{{\rm min},i}$ (can be negative) and $u_{{\rm max},i}$ are the minimum and maximum of the $i$-th element in the labels of all manifold data, respectively, $\lambda_1$ and $\lambda_2$ are hyper-parameters. We set $\lambda_1=0.5$ and $\lambda_2=2$ in our study.

Then for each data $\boldsymbol{x}(t_0):=\boldsymbol{x}(t)$ after filteration, given the evolution time sequence $\boldsymbol{\tau} = [\tau_1, \tau_2, \cdots, \tau_k]$, we simulate a short reaction trajectory and sample several data points from the trajectory. We choose $k+1$ samples in each trajectory, namely $\boldsymbol{x}(t_i):=\boldsymbol{x}(t_{i-1}+\tau_i)$, $i=0, 1, \cdots, k$. The time sequence for chemical kinetics model is set as $\boldsymbol{\tau}=[10^{-6}, 10^{-6}, 10^{-5}, 10^{-5}, 10^{-4}, 10^{-4}, 10^{-3}, 2\times10^{-3}, 5\times 10^{-3}]$ s.
All the data after temporal evolution will be merged together. Then the new dataset is filtered again by the manifold data $\mathbb{D}_{\rm Init}$. Finally, we obtain the training dataset $\mathbb{D}_{\rm EMCS}$. During the training process, the DNN hidden neurons for chemical kinetics is set as 3200, 1600, 800, 400.

\subsection{Computational Singular Perturbation}
\label{sec:csp}

Computational singular perturbation (CSP) is a method for analyzing dynamical systems. CSP theory is based on the assumption that typical multi-scale dynamics system could be regarded as the combination of fast and slow modes. Fast modes limit the time scale of numerical integration which yielding the stiffness, while the slow modes allow larger time scales that could be treated with explicit schemes. Then the CSP method aims to decouple the slow time modes and the fast time modes. The CSP framework is as follows.

Consider the autonomous stiff ODE system in the form:

\begin{align}
    \frac{\mathrm{d}\boldsymbol{x}}{\mathrm{d}t} = \boldsymbol{\omega}(\boldsymbol{x})
\end{align}
with the initial value $\boldsymbol{x}_0$ and the state vector $\boldsymbol{x}\in \mathbb{R}^{N}$. Let matrix $A=[\boldsymbol{a}_1, \boldsymbol{a}_2, \cdots, \boldsymbol{a}_N]$, the column vectors $\boldsymbol{a}_i\in\mathbb{R}^{N}$, and $B=[\boldsymbol{b}_1, \boldsymbol{b}_2, \cdots, \boldsymbol{b}_N]$ is the inverse matrix. The row vectors $\boldsymbol{b}_i$ satisfy the orthogonality property 
\begin{align}
    \boldsymbol{b}_i\cdot\boldsymbol{a}_j = \delta_{ij}
\end{align}
We rewrite the source term of ODE $\boldsymbol{\omega}$:
\begin{align}
    \boldsymbol{\omega} = I\cdot\boldsymbol{\omega} = AB \cdot \boldsymbol{\omega} = A\cdot\boldsymbol
{f} = \sum_{i=1}^{N} f_i \boldsymbol{a}_i
\end{align}
where $\boldsymbol{f}=B\boldsymbol{\omega}$ and $f_i$ is the signed amplitude of the projection of $\boldsymbol{\omega}$ on the basis vector $\boldsymbol{b}_i$ since
\begin{align}
    f_i = b_i\cdot\boldsymbol{\omega}
\end{align}
CSP method seeks a set of appropriate vectors $\boldsymbol{a}_i$ such that the source term can be spanned by these vectors which are referred to as CSP basis vectors. Estimating the CSP basis vectors is commonly achieved by utilizing the right eigenvectors of the Jacobian matrix $J = \frac{\partial\boldsymbol{\omega}}{\partial \boldsymbol{x}}$. Then we have the  eigen-decomposition property 
\begin{align}
   JA = A\Lambda, J = A\Lambda B
\end{align}
where $\Lambda=diag(\lambda_1, \lambda_2, \cdots, \lambda_N)$ is the diagonal matrix of the eigenvalues of $J$.
Further, we can study the dynamics of $\boldsymbol{f}$
\begin{align}
   \frac{\mathrm{d}\boldsymbol{f}}{\mathrm{d}t} = \frac{\mathrm{d}B } {\mathrm{d}t} \boldsymbol{\omega} +  B \frac{\mathrm{d}\boldsymbol{\omega}}{\mathrm{d}t} = BJ\boldsymbol{\omega} = \Lambda\cdot \boldsymbol{f}  
\end{align}
where the $\mathrm{d}B/\mathrm{d}t$ is set to zero. Each $f_i$ satisfies the ODE form:
\begin{align}
   \frac{\mathrm{d}f_i} {\mathrm{d}t} = \lambda_i f_i, \quad i=0,1,\cdots, N  
   \label{eq:fi_ODE}
\end{align}
The reciprocal of eigenvalues $\tau_i:=1/|\lambda_i|$ show the approximation timescales of the dynamic system. 
We arrange the time scales $\tau_i$ in ascending order,
\begin{align}
   \tau_1 \leq \tau_2\leq \cdots\tau_{M}\ll\tau_{M+1}\leq\cdots\leq\tau_{N}
\label{eq:tau_seq}
\end{align}
Hence mode 1 associated with $\tau_1$ is the fastest and model N with $\tau_N$ is slowest. A very small time scale $\tau_i$ implies a large eigenvalue $\lambda_i$, and the corresponding amplitude $f_i$ will decay exponentially if the real part of the eigenvalue is negative. So we can split $\boldsymbol{\omega}$ into slow part and fast part,
\begin{align}
   \boldsymbol{\omega} = \underbrace{\sum_{i=1}^{M} f_i \boldsymbol{a}_i}_{\boldsymbol{\omega}_{\text{fast}}\approx 0     } + \underbrace{\sum_{i=M+1}^{N} f_i \boldsymbol{a}_i} _{\boldsymbol{\omega}_{\text{slow}}\approx 0}
\end{align}
and $f_{M+1}$ associated with $\tau_{M+1}$ mentioned in Eq. \ref{eq:tau_seq} is the fastest one within slow modes. To decouple fast modes and slow modes, we need to identify $\tau_{M+1}$. The motivation is that the integration of $\boldsymbol{\omega}_{\text{fast}}$ over time scale $\tau_{M+1}$ is negligible, namely,
\begin{align}
    \delta x^{j}_{\text{fast}} = \left\vert \int^{\tau_{M+1}}_{0} \omega_{\text{fast}}^{j} \mathrm{d} t \right\vert = \left|\sum_{i=1}^{{M}} a_i^j f_i \frac{e^{\lambda_i \tau_{M+1}} -1}{\lambda_i}\right|<\delta x_{\mathrm{err}}^j, \quad j=1, \cdots, N
\end{align}
where $\delta x_{\mathrm{err}}^j$ is the $j$-th component of $\delta \boldsymbol{x}_{\mathrm{err}}=\text{tol}_{\mathrm{rel}}|\boldsymbol{x}|+\text { tol}_{\text{abs}}$ and $\boldsymbol{x}$ is the local state vector. In our study, $\tau_{csp} = \tau_{M+1}$ and the hyper-parameters $\text{tol}_{\text{rel}} = 10^{-4}$, $\text{tol}_{\text{abs}} = 10^{-10}$ for chemical kinetics model.

\section{Additional Results}

\subsection{Details about Electronic Dynamical Process}
\label{sec:modulator_appd}

The electronic dynamical process (ring modulator), originating from a electrical circuit, is chosen as a test problem for comparing the integration methods in previous study \cite{Timofeev2022modulator}. It contains 15 composition elements with significant nonlinearity and significant stiffness. Hence, we test our DeePODE method via this problem. 

The dynamical system receives a low-frequency signal $U_{in1}$ and a high frequency signal $U_{in2}$, then produce generates a mixed
signal at the output. The state vector can be expressed by voltages $U_i$ and currents $I_i$:
\begin{align}
    y=\left(U_1, U_2, U_3, U_4, U_5, U_6, U_7, I_1, I_2, I_3, I_4, I_5, I_6, I_7, I_8\right)
\end{align}
The governing equation for this non-autonomous dynamical system is as follows.
\begin{equation}
\begin{gathered}
y_1^{\prime}=\frac{1}{C}\left(y_8-0.5 y_{10}+0.5 y_{11}+y_{14}-\frac{1}{R} y_1\right), \\
y_2^{\prime}=\frac{1}{C}\left(y_9-0.5 y_{12}+0.5 y_{13}+y_{15}-\frac{1}{R} y_2\right), \\
y_3^{\prime}=\frac{1}{C_s}\left(y_{10}-q\left(U_{D 1}\right)+q\left(U_{D 4}\right)\right), \\
y_4^{\prime}=\frac{1}{C_s}\left(-y_{11}+q\left(U_{D 2}\right)-q\left(U_{D 3}\right)\right), \\
y_5^{\prime}=\frac{1}{C_s}\left(y_{12}+q\left(U_{D 1}\right)-q\left(U_{D 3}\right)\right), \\
y_6^{\prime}=\frac{1}{C_s}\left(-y_{13}-q\left(U_{D 2}\right)+q\left(U_{D 4}\right)\right), \\
y_7^{\prime}=\frac{1}{C_p}\left(-\frac{1}{R_p} y_7+q\left(U_{D 1}\right)+q\left(U_{D 2}\right)-q\left(U_{D 3}\right) -q\left(U_{D 4}\right)\right), \\
y_8^{\prime}=-\frac{1}{L_h} y_1 \\
y_9^{\prime}=-\frac{1}{L_h} y_1, \\
y_{10}^{\prime}=\frac{1}{L_{s 2}}\left(0.5 y_1-y_3-R_{g 2} y_{10}\right), \\
y_{11}^{\prime}=\frac{1}{L_{s 3}}\left(-0.5 y_1+y_4-R_{g 3} y_{11}\right), \\
y_{12}^{\prime}=\frac{1}{L_{s 2}}\left(0.5 y_2-y_5-R_{g 2} y_{12}\right), \\
y_{13}^{\prime}=\frac{1}{L_{s 3}}\left(-0.5 y_2+y_6-R_{g 3} y_{13}\right), \\
y_{14}^{\prime}=\frac{1}{L_{s 1}}\left(-y_1+U_{i n 1}(t)-\left(R_i+R_{g 1}\right) y_{14}\right), \\
y_{15}^{\prime}=\frac{1}{L_{s 1}}\left(-y_2 - \left(R_c+R_{g 1}\right) y_{15}\right), 
\end{gathered} 
\end{equation}
The initial value is set as 
\begin{align}
    y_i(0)=0, ~~i =1, \cdots, 15, ~~~0 \leq t \leq 10^{-3} s
\end{align}
Functions $U_{D 1}, U_{D 2}, U_{D 3}$ , $U_{D 4}, q, U_{in1 }$ , $U_{i n 2}$ are defined by equations
\begin{equation}
\begin{gathered}
U_{D 1}=y_3-y_5-y_7-U_{i n 2}(t), \\
U_{D 2}=-y_4+y_6-y_7-U_{i n 2}(t), \\
U_{D 3}=y_4+y_5+y_7+U_{i n 2}(t), \\
U_{D 4}=-y_3-y_6+y_7+U_{i n 2}(t), \\
q(U)=\gamma\left(e^{\delta U}-1\right), \\
U_{i n 1}(t)=0.5 \sin (2000 \pi t), \\
U_{i n 2}(t)=2 \sin (20000 \pi t) .
\end{gathered}
\end{equation}
The parameters mentioned above are:
$$
\begin{gathered}
C=1.6 \cdot 10^{-8}, \quad C_s=2 \cdot 10^{-12}, C_p=10^{-8} \\
L_h=4.45, L_{s 1}=0.002, L_{s 2}=5 \cdot 10^{-4} \\
L_{s 3}=5 \cdot 10^{-4}, \gamma=40.67286402 \cdot 10^{-9} \\
R=25000, R_p=50, R_{g 1}=36.3, R_{g 2}=17.3 \\
R_{g3} = 17.3, R_i = 50, R_c = 600, \delta = 17.7493332.
\end{gathered}
$$

\begin{figure}[htbp]
    \centering
    \includegraphics[width=1\linewidth]{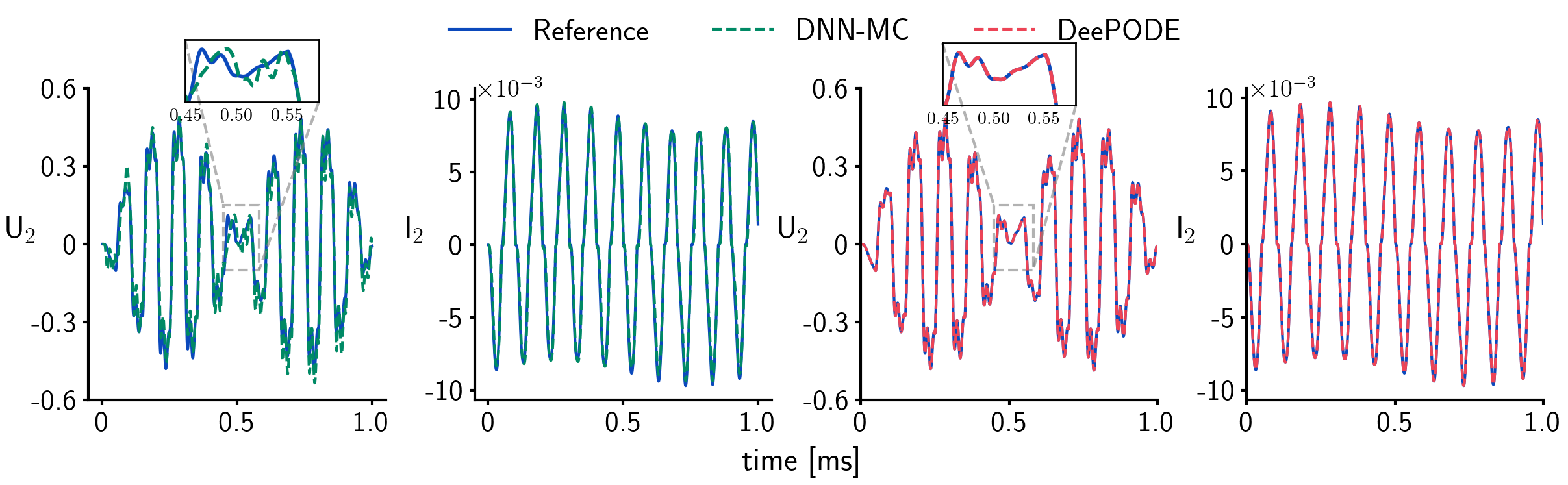}
    \caption{\textbf{Ring modulator model}. The DNN-MC/DeePODE predictions and direct integration results for $U_2$ and $I_2$ are evaluated with an initial condition of $y_i(0) = 0$ over the time interval $0 \leq t \leq 10^{-3}~\rm{s}$. The numerical integration is performed using the Runge-Kutta scheme. } 
    \label{fig:modulator_appd}
\end{figure}

In the EMCS stage, we first define the initial condition range as $U_i, I_i \in [-10^{-3}, 10^{-3}]$ for manifold sampling, and randomly select 800 initial state vectors. Each state is then simulated for $1~\rm{ms}$, with 1,000 data points collected at a time interval of $10^{-6}~\rm{s}$. All the data points are combined and referred to as $\mathbb{D}_{\text{MF}}$.
For the MC sampling, the range is set as $I_1 \in [-10^{-5}, 10^{-7}]$, $I_2 \in [-10^{-6}, 10^{-6}]$, with the remaining components in $[-10^{-3}, 10^{-3}]$. Subsequently, 5,000 points are randomly sampled and simulated for $k = 30$ steps. Each evolution time $\tau_j, 0 \leq j \leq k$, is randomly selected from $\tau_j \in [0, 0.02]~\rm{s}$.
Finally, all data are filtered based on the range of temporal change $[\lambda_1 u^{\min}_i, \lambda_2 u^{\max}_i]$ over $\mathbb{D}_{\text{MF}}$, where $\lambda_1 = \lambda_2 = 10$. Ultimately, we obtain the training dataset $\mathbb{D}_{\text{EMCS}}$.

We further present the comparison between DNN-MC/DeePODE predictions and direct integration results for $U_2$ and $I_2$ (Fig.~ \ref{fig:modulator_appd}). Our findings reveal that, in the high-frequency oscillation range of the ring modulator, DNN-MC faces difficulties in making accurate predictions, whereas DeePODE consistently provides accurate results. This demonstrates that DeePODE effectively captures the multiscale dynamics inherent in the system.

\subsection{Reaction-Diffusion Systems}
\label{sec:combustion_appd}

\subsubsection{Governing Equation}

Assuming ideal gases and perfect mixtures, the conservation equations of mass, momentum, species, and energy of \textit{EBI-DNS} \cite{zirwes2020quasi} used in this work, are given as follows:
\begin{equation}
\frac{\partial \rho}{\partial t} + \nabla \cdot (\rho \Vec{u}) = 0,
\end{equation}
\begin{equation}
    \frac{\partial (\rho \Vec{u})}{\partial t} + \nabla \cdot (\rho \Vec{u} \Vec{u}) = - \nabla p + \nabla \cdot \Vec{\tau} + \rho \Vec{g}, \label{eq:momentum}
\end{equation}
\begin{equation}
    \frac{\partial (\rho Y_k)}{\partial t} + \nabla \cdot (\rho(\Vec{u}+\Vec{u_c})Y_k) = \dot{\omega}_k - \nabla \cdot \Vec{j}_k, \quad k=1 \ldots N-1, \label{eq:species}
\end{equation}
\begin{equation}
    \frac{\partial(\rho(h_s + \frac{1}{2}\Vec{u} \cdot \Vec{u}))}{\partial t} + \nabla \cdot (\rho \Vec{u}(h_s + \frac{1}{2}\Vec{u}\cdot\Vec{u})) = - \nabla \cdot \Vec{\dot{q}} + \frac{\partial p}{\partial t} - \sum_{k=1}^{N} \mathring{h}_k \dot{\omega}_k, \label{eq:energy}
\end{equation}
where $t$ is time, $\rho$ is the fluid density, $\Vec{u}$ is the fluid velocity. In Eq.~[\ref{eq:momentum}], $p$ is fluid pressure, $\Vec{g}$ is the gravitational acceleration and $\Vec{\tau}$ is the viscous stress tensor of Newtonian fluid using Stokes assumption:
\begin{equation}
    \Vec{\tau} = \mu (\nabla \Vec{u} + (\nabla \Vec{u})^{T} - \frac{2}{3}\Vec{I}\nabla \cdot \Vec{u}),
\end{equation}
$\Vec{I}$ is the identity tensor, $\mu$ is the dynamical viscosity, which is calculated from the Chapman-Enskog solution of kinetic gas theory. In Eq.~[\ref{eq:species}], $N$ is the number of species.  $\Vec{u}_c$ is the correction velocity that forces the sum of all diffusive fluxes $\Vec{j}_k$ to be zeros: $\Vec{u}_c = - \frac{1}{\rho} \sum_{k=1}^{N} \Vec{j}_k$, where $\Vec{j}_k$ is the diffusive mass flux of $k$-th species. $Y_k$ is the mass fraction of $k$-th species and $\dot{w}_k$ is its reaction rate and obtained by
\begin{equation}\label{eq:w_k}
    \dot{w}_k = M_k\sum_r (\nu_{k,r}^{''}-\nu_{k,r}^{'})(k_r^{'} \prod_k C_k^{\nu_{k,r}^{'}}-k_r^{''}\prod_k C_k^{\nu_{k,r}^{''}})
\end{equation}
in traditional solver where $C_k=\frac{Y_k}{M_k}\rho$ is the molar  concentration of species $k$, $\nu_{k,r}^{'}(\nu_{k,r}^{''})$ is the forward (reverse) stoichiometric coefficient of species $k$ and reaction $r$, $k_r^{'}(k_r^{''})$ is the forward (reverse) rate constant of reaction $r$ which is computed by Arrhenius equation. \emph{In DeePODE, Eq.~[\ref{eq:w_k}] will be computed by DNN model}. In \textit{EBI-DNS}, there are two diffusion models, the full multi-component diffusion including Soret diffusion and the mixture-averaged diffusion model (Hirschfelder-Curtiss approximation). For all \textit{EBI-DNS} examples in this paper, we compute the diffusion mass flux via a mixture-averaged diffusion model:
\begin{equation}
    \Vec{j}_k = -\rho D_{m, k}^{mole} \nabla Y_k - Y_k \rho D_{m,k}^{mole} \frac{1}{\Bar{M}} \nabla \Bar{M},
\end{equation}
\begin{equation}
    \quad  D_{m, k}^{mole} = \frac{1-Y_k}{\sum_{j \neq k}\frac{X_j}{D_{j, k}}},
\end{equation}
where $M_k$ is the molar mass of k-th species, $\Bar{M}$ is the mean molecular weight. ${D_{j,k}}$ is the binary diffusion coefficient between $k$-th species and $j$-th species, which is calculated from the Chapman-Enskog solution of kinetic gas theory.

In Eq.~[\ref{eq:energy}], the transport of energy is formulated in terms of the total sensible enthalpy $h_s + \frac{1}{2}\Vec{u} \cdot \Vec{u}$, $h_s$ is the sensible enthalpy of mixture for ideal gases, $h_s = \sum_{k=1}^{N}Y_k h_{s,k}$, where $h_{s,k}=h_k - \mathring{h}_k$, $\mathring{h}_k=h_k(298 {\rm K})$ is the enthalpy of formation of the $k$-th species. The divergence of the negative energy flux $-\nabla \cdot \dot{\Vec{q}}$:
\begin{equation}
    -\nabla \cdot \dot{\Vec{q}} = \nabla \cdot (\frac{\lambda}{c_p}\nabla h_s) - \sum_{k=1}^{N}\nabla \cdot (\frac{\lambda}{c_p}h_{s,k}\nabla Y_k) - \sum_{k=1}^{N}\nabla \cdot (h_{s,k}\hat{\Vec{j}}_k), \label{eq:energyflux}
\end{equation}
where $\lambda$ is the heat conductivity of the mixture, which is calculated from the Chapman-Enskog solution of kinetic gas theory. $c_p$ is the isobaric heat capacity and the corrected diffusive mass flux $\hat{\vec{j}}_k = \vec{j}_k - Y_k\sum_{i=i}^{N}\Vec{j}_i$.
Notice that under the condition that the gases are ideal gases and all species have the same temperature, by the Fourier second law, we have
\begin{equation}
    \nabla \cdot (\lambda \nabla T) =  \nabla \cdot (\frac{\lambda}{c_p}\nabla h_s) - \sum_{k=1}^{N}\nabla \cdot (\frac{\lambda}{c_p}h_{s,k}\nabla Y_k),
\end{equation}
thus, Eq.~[\ref{eq:energyflux}] can be written as 
\begin{equation}
    -\nabla \cdot \dot{\Vec{q}} = \nabla \cdot (\lambda \nabla T) - \sum_{k=1}^{N}\nabla \cdot (h_{s,k}\hat{\Vec{j}}_k).
\end{equation}

\subsubsection{Chemical Mechanisms}
\paragraph{GRI-Mech 3.0 methane mechanism.}
GRI-Mech 3.0, created by the Berkeley Combustion Team, is a detailed kinetic mechanism for methane/air combustion including nitric species for NOx predictions with 53 species and 325 reactions~\footnote{Gregory P. Smith, David M. Golden, Michael Frenklach, Nigel W. Moriarty, Boris Eiteneer, Mikhail Goldenberg, C. Thomas Bowman, Ronald K. Hanson, Soonho Song, William C. Gardiner, Jr., Vitali V. Lissianski, and Zhiwei Qin, \url{http://combustion.berkeley.edu/gri-mech/version30/text30.html}.}. It is challenging to obtain a surrogate model of GRI-Mech 3.0 because of the large number of species. For example, \cite{chi2021fly} used an on-the-fly scheme to overcome this difficulty but with large computation cost, i.e.,  frequently retraining the neural network with data from the current simulation. 

\paragraph{DRM19 methane mechanism.}
DRM19 is a reduced reactions subset of full GRI-Mech 1.2, including 21 species and 84 reactions~\footnote{A. Kazakov, M. Frenklach, Reduced Reaction Sets based on GRI-Mech 1.2, a 19-species reaction set, \url{http://www.me.berkeley.edu/drm/}.}. The reduced mechanism has been compared with GRI-Mech 1.2 in ignition delay and laminar flame speed simulations, demonstrating reasonable accuracy and reduced computational costs.



\paragraph{n-Heptane mechanism.} An in-house reduced n-heptane mechanism consisting of 34 species and 191 reactions \footnote{The reduced mechanism can be found at \url{https://github.com/intelligent-algorithm-team/intelligent-combustion}.} from  a detailed mechanism with 116 species and 830 reactions \cite{chaos2007high}. The reduced mechanism obtained by the method DeePMR \cite{DeePMR2022} is based on the original detailed n-heptane with 116 species and 830 reactions \cite{chaos2007high}.

\subsubsection{CFD Codes}

To demonstrate the good usability and portability of DeePODE, we integrated the DeePODE model into different CFD codes, including both C++ and FORTRAN implementations, to perform combustion simulations.

\paragraph{\textit{EBI-DNS}.} 

\textit{EBI-DNS} \cite{zirwes2020quasi} is a C++ package implemented within the OpenFOAM-5 framework, designed for the accurate and efficient simulation of compressible reactive flows. It integrates the extended $\rm{rhoReactingBuoyantFoam}$ solver with the open-source library Cantera to compute thermodynamic and transport properties. The ODE integration is handled using the Sundials CVODE solver, ensuring robust and efficient numerical solutions.

\paragraph{ASURF.}

\textit{ASURF} \cite{chen2009effects,zhang2021studies} is an in-house code developed in FORTRAN for the \textbf{A}daptive \textbf{S}imulation of \textbf{U}nsteady \textbf{R}eacting \textbf{F}low. It employs the finite volume method (FVM) to solve the conservation equations for multi-species mixtures. Thermo-physical and transport properties are computed using the CHEMKIN and TRANSPORT libraries, while chemical kinetics are integrated using the built-in VODE ODE solver.

\paragraph{DeepFlame.}
\textit{DeepFlame} \cite{mao2022deepflame} is a deep learning empowered computational fluid dynamics package for single or multiphase, laminar or turbulent, reacting flows at all speeds. It aims to provide an open-source platform to combine the individual strengths of OpenFOAM, Cantera and PyTorch libraries for deep learning assisted reacting flow simulations.

\subsubsection{More Examples with \textit{EBI-DNS} Code}

\paragraph{0D auto-ignition.}

\begin{figure*}[htp]
    \centering
    \includegraphics[width=1\textwidth]{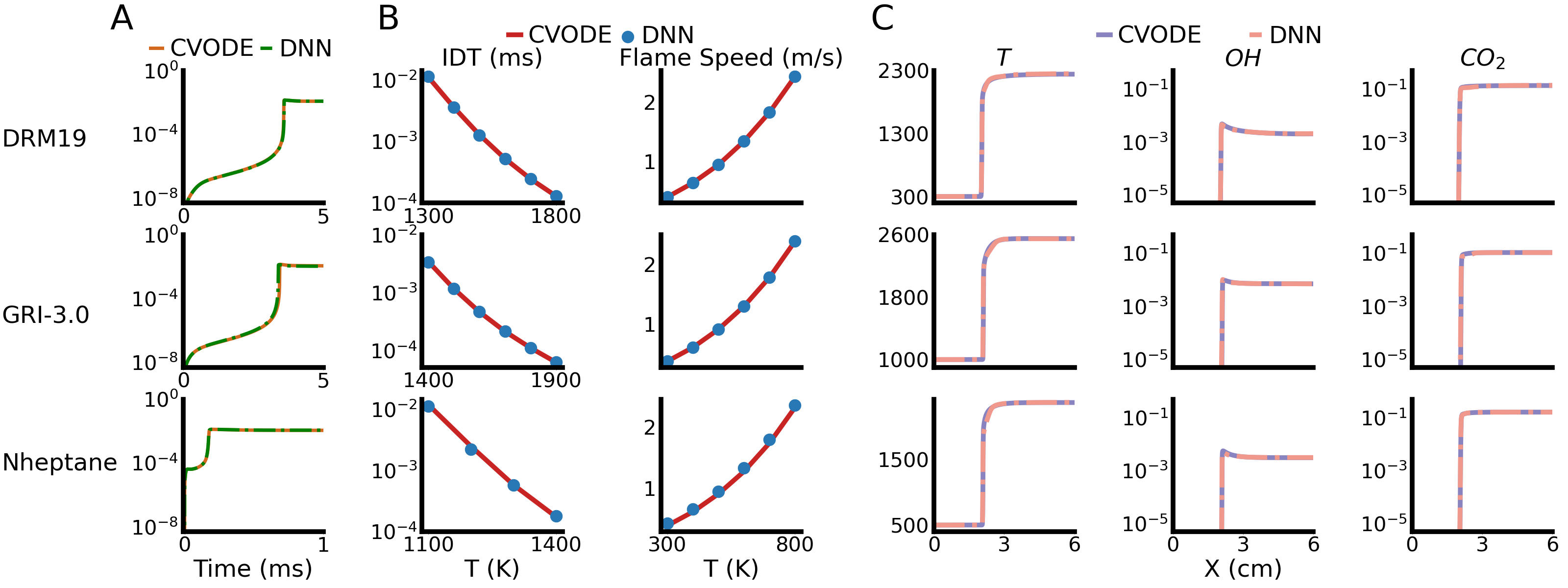}
    \caption{ \textbf{Zero-dimensional and One-dimensional examples of chemical kinetic system}. (A) Trajectory of OH with initial temperature 1400K, pressure 1 atm, and equivalence ratio 1. (B) A comparison of ignition delay time and flame speeds under different initial temperatures (abscissa) and equivalence ratio 1 at $1 {\rm atm}$. Blue points are computed by iteratively using DeePODE model and red curves by CVODE.  (C) Temperature, ${\rm OH}$, ${\rm CO}_2$ over the one-dimensional space at initial temperature 300K (DRM19), 1000K (GRI) and 500K (n-heptane).}\label{fig:0d1d}
\end{figure*}

The basic test example is a zero-dimensional constant-pressure auto-ignition one with given initial temperature, pressure, and equivalence ratio. The DeePODE model is used iteratively to obtain a complete combustion trajectory and compared with results from Cantera. For illustration, Fig. \ref{fig:0d1d}(A) shows the OH trajectory, where the ones obtained by the DeePODE models nearly overlap with those obtained by Cantera over all orders of magnitude. For different initial conditions, as shown in the first column in Fig. \ref{fig:0d1d}(B), the ignition delay times (the time point where the temperature change rate reaches a maximum) under different initial conditions for the DeePODE models and Cantera are in excellent agreement.

\paragraph{1D freely propagating premixed flame.}

The one-dimensional example is a freely propagating premixed flame. An inlet on the left provides fresh gases in air at 1 atm and an equivalence ratio of one to the domain, which exits at the outlet on the right. The left half of the domain is initially filled with the unburnt mixture, and the right half is filled with the completely burnt mixture. After some time, the flame reaches a stationary state. The physical length is $0.06 {\rm m}$ m with 4335 cells. Statistically, we compare the flame speed (measured rate of expansion of the flame front in a combustion reaction) under different initial temperatures at an equivalence ratio of 1 and 1 atm. As shown in the second column in Fig. \ref{fig:0d1d}(B),  the DeePODE model can accurately predict flame speeds in all test examples. As shown in Fig. \ref{fig:0d1d}(C), the distributions of temperature and mass fractions of OH, CO$_2$, and CH$_4$ over space are very consistent between results obtained by the DeePODE model and CVODE in Cantera.

\paragraph{2D counter flow.}

\begin{figure*}[htp]
    \centering
    \includegraphics[width=1\textwidth]{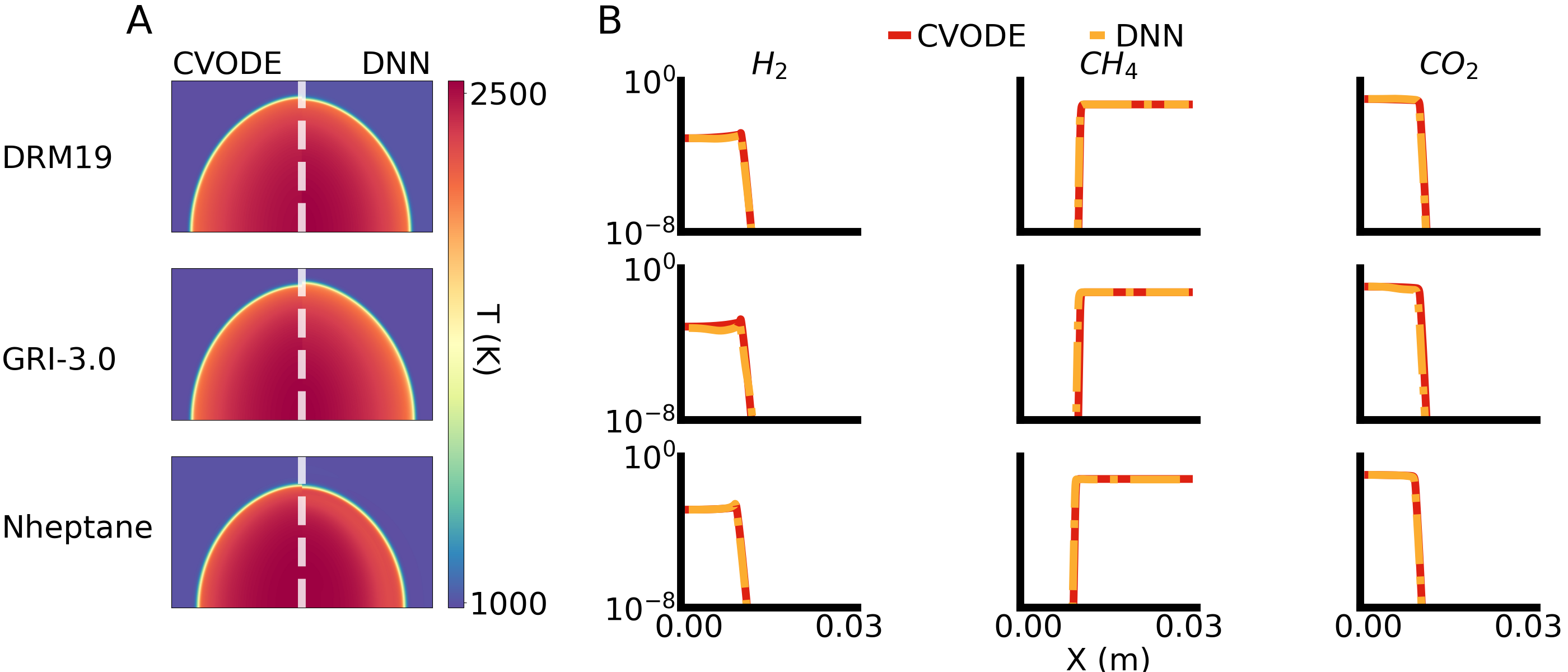}
    \caption{\textbf{Two-dimensional counter flow}. (A) shows the temperature distribution by \textit{EBI-DNS} with CVODE and DNN, respectively, at $1 {\rm ms}$ of counter flow. (B) is the comparison between CVODE and DNN at one-dimensional cross section of counter flow.}\label{fig:counterflow}
\end{figure*}

In the premixed counter flow, the computation domain size is $2 {\rm cm}\times 3{\rm cm}$ with 960,000 cells. Premixed gas with $\phi = 1$ is injected from both sides into the middle.
We set the ignition area as a circle with a radius of $0.02{\rm cm}$. After some time, due to the constant introduction of premixed gas to the left and right sides, the flame no longer spreads to the two sides but only to the front. As shown in Fig.~\ref{fig:counterflow}(A), we compare the results at $1 {\rm ms}$ and find that the DNN models can accurately capture the flame structure and propagation for DRM19, GRI 3.0, and n-heptane. For illustration, Fig.~\ref{fig:counterflow}(B) shows the distribution of temperature, mass fractions of ${\rm H}_2$, ${\rm CH}_4$, ${\rm CO}_2$ at a one-dimensional cross-section.

\begin{figure}[htp]
    \centering
    \includegraphics[width=1\linewidth]{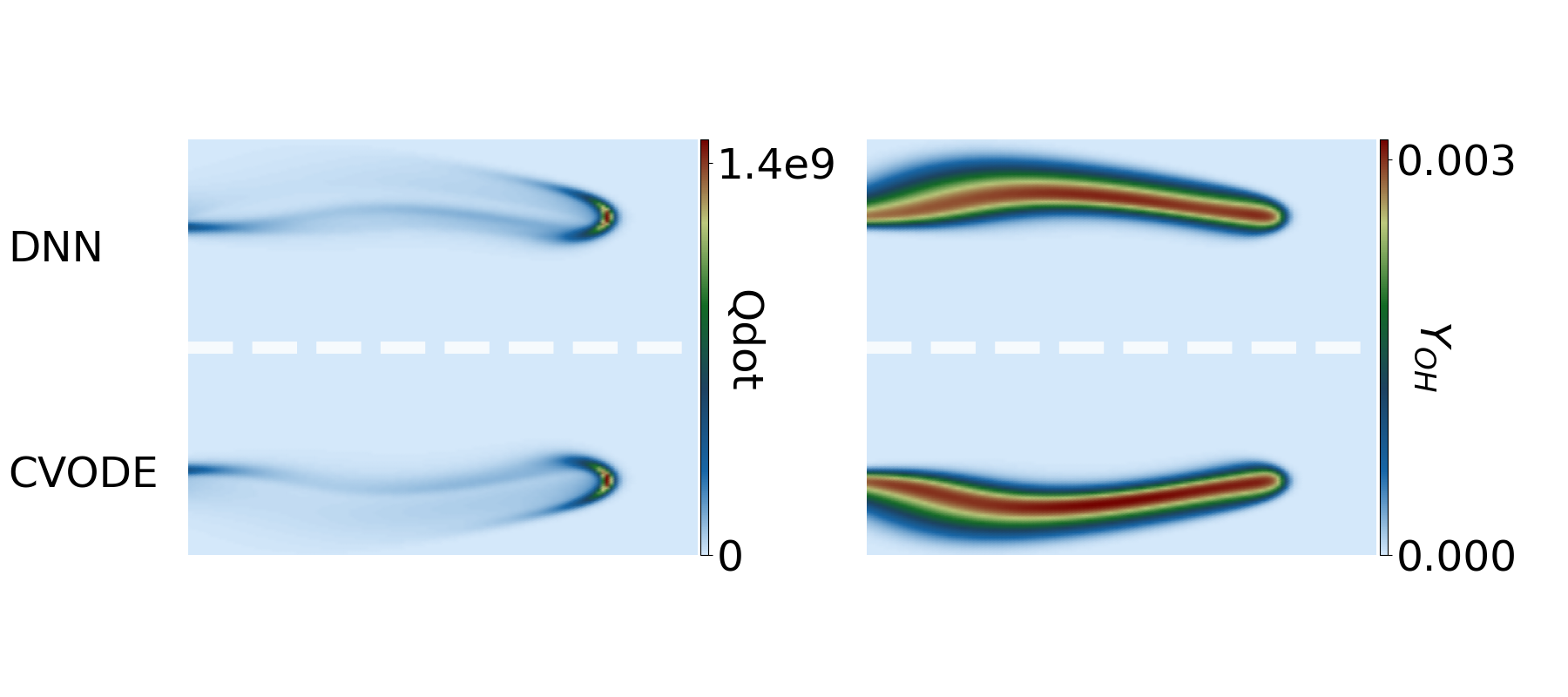}
    \caption{\textbf{Triple flame}. The results (heat release rate and OH distribution) by  \textit{EBI-DNS} with CVODE and DeePODE (DNN), respectively, at $6 {\rm  ms}$  for DRM19.}\label{fig:triple}
\end{figure}

\paragraph{2D triple flame.}
In the triple flame case, it consists of a $3{\rm cm}\times5 {\rm cm}$ computation domain with a jet surrounded by an air co-flow, which contains 150,000 cells. We set a high-temperature domain in the central of the computation domain to ignite.The specific condition can be seen in the picture. As shown in Fig. \ref{fig:triple}, for both the heat release rate and OH distribution at $6 {\rm  ms}$  for DRM19, results by  CVODE and DNN agrees very well. 
 
\subsubsection{Example with \textit{ASURF} Code}
\paragraph{One-dimensional freely propagating premixed flame.} 
This example is similar to the one in Fig. \ref{fig:0d1d} with GRI, except that the flow comes in from the right hand side. Fig. \ref{fig:1d-ASURF}(A)
shows the temporal evolutions of flame front location and propagation speed at $T= 1600 K$, $P = 1 atm$, $\phi = 1.1$ using DNN or VODE as the integrator. As shown in Fig. \ref{fig:1d-ASURF}(A), these results of DNN and VODE agree very well. 

\subsubsection{Example with \textit{DeepFlame} Code}
\begin{figure}[htbp]
	\centering
        \includegraphics[width=1\linewidth]{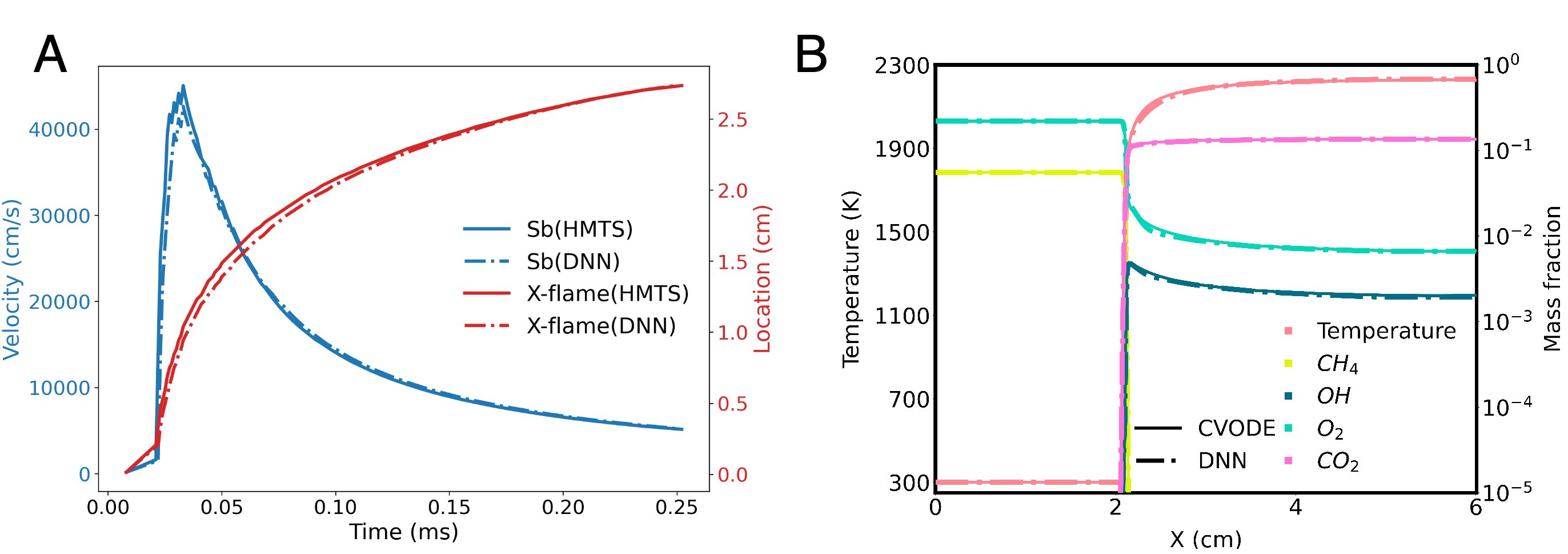}
 
  \caption{\textbf{One-dimensional freely propagating premixed flame}. The figure shows flame front location and propagation speed at $T= 1600 {\rm K}$, $P = 1 {\rm atm}$, $\phi = 1$ using DNN or for (A) VODE as the integrator with ASURF code, and for (B) CVODE as the integrator with \textit{DeepFlame} code.} \label{fig:1d-ASURF}
\end{figure}

\paragraph{One-dimensional freely propagating premixed flame.} 
This example is same to the one in Fig. \ref{fig:0d1d} with DRM19. Fig. \ref{fig:1d-ASURF}(B)  
shows the distribution of temperature and mass fraction of ${\rm O}_2$, ${\rm CO}_2$, OH over space at $T= 300 K$, $P = 1 atm$, $\phi = 1$ using DNN or CVODE as the integrator. As shown in Fig. \ref{fig:1d-ASURF}(B), these results of DNN and CVODE agree very well.

\section{Additional Analysis}

\subsection{Solver Indicator}

We define a distance function of a data $y$ and a dataset $X$ to evaluate the feasibility of using a DNN model. Since the distance between the data and the training set is difficult to obtain directly due to its non-convexity and high dimensionality, we use the following probability density distribution to approximate the distance of the data to the training set.
 \begin{equation}\label{eq:p}
     p(y,X) = \prod_{i=1}^N \sqrt{2\pi}\phi\left(\frac{y_i-\mu(X) } {\sigma(X)} \right) = \prod_{i=1}^N e^{-\frac{1}{2}\left(\frac{y_i-\mu(X)}{\sigma{X}}\right)^2}
 \end{equation}
 where $\phi$ denotes the probability density function of the standard normal distribution, and $\mu(X)$,$\sigma(X)$ is the 
mean and standard deviation, respectively. The smaller is $p(y,X)$, the further distance is data from the training set. The confidence interval is dependent on the DNN model. 

\begin{figure}[htbp]
    \centering
    \includegraphics[width=0.7\linewidth]{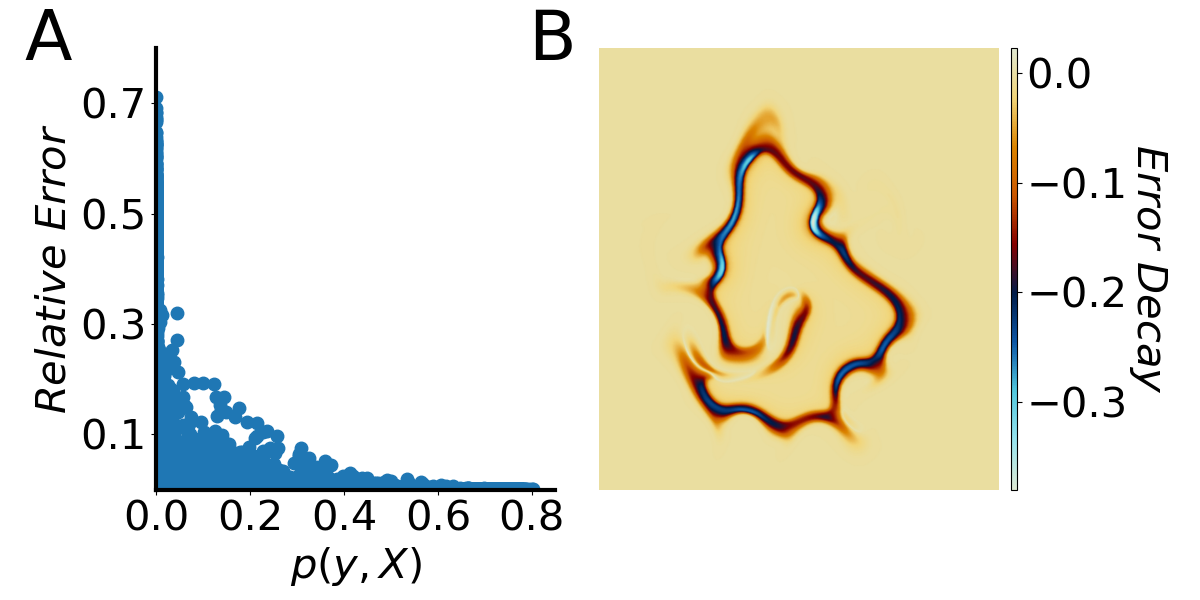}
    \caption{ \textbf{(A)} The distribution of temperature's relative error against to $d(y,X)$, contains 500000 data points from random sampling. \textbf{(B)} The comparison between whether implementing solver indicator on turbulent ignition example. The color means the decrease of relative error brought by the indicator.  }\label{fig:pos_error}
\end{figure}

We implement the solver indicator on the DRM19-DeePODE model to verify its effect.  
 Fig. \ref{fig:pos_error}(A) illustrates that most data points with large errors are located in domains farther away from the training set. Although the assumption of normal distribution is not justified, it empirically turns out to be a good indicator for the distance of the data to the training set. We choose the confidence integral $[0.5,1]$ and test it on the turbulent ignition example. Fig. \ref{fig:pos_error}(B) states that the indicator can effectively decrease the errors. Subsequent works can have hybrid methods with both traditional methods and DeePODE based on this indicator to improve accuracy.

\subsection{Additional Characteristic Time Analysis}

We also investigate the characteristic time distribution in both a two-dimensional predator-prey system and a 15-dimensional modulator case.

\begin{figure}[htbp]
	\centering
        \includegraphics[width=1\textwidth]{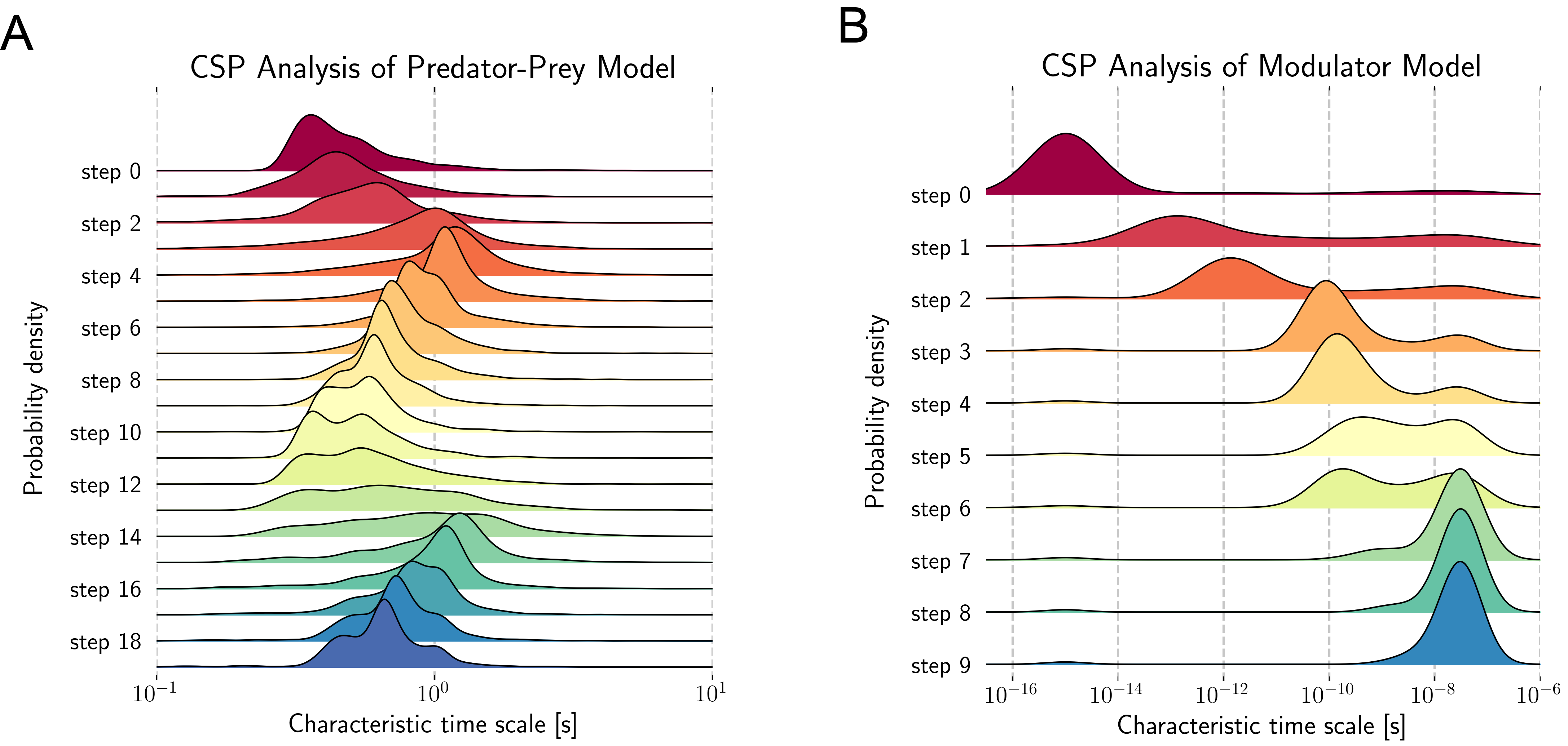}
\caption{\textbf{
EMCS characteristic time scales distribution}. \textbf{(A)} Comparison of the characteristic time scales distribution at different evolution step in EMCS for 2D predator-prey model. \textbf{(B)}  Comparison of the characteristic time scales distribution at different evolution step in EMCS for 15-dimensional ring modulator model.  
} 
\label{fig:csp_time_preadtor_modulator}
\end{figure}

For the predator-prey dynamics, we randomly select 2,000 samples within the range $x_1 \in [0, 5]$, $x_2 \in [0, 5]$. The evolution time sequence for the EMCS method is set as $\boldsymbol{\tau} = [\tau_i], i = 0, \dots, 19$, with each $\tau_i = 0.5~\rm{s}$. Since the Jacobian matrix has only two eigenvalues, we estimate the characteristic timescale as $1/\lambda_{\max}$ for state $\boldsymbol{x}_t$, where $\lambda_{\max}$ represents the largest eigenvalue of the local Jacobian matrix. As shown in Fig. \ref{fig:csp_time_preadtor_modulator}(A), we visualize the characteristic timescale distribution from step 0 to step 19 in the EMCS method. We found that EMCS effectively captures multiple characteristic timescales with only a few evolution steps. Furthermore, the characteristic time distributions exhibit periodic behavior as the steps progress, which corresponds to the inherently periodic nature of the predator-prey solution.

For the ring modulator model, the range for the MC method is defined as $I_1 \in [-10^{-5}, 10^{-7}]$, $I_2 \in [-10^{-6}, 10^{-6}]$, with the remaining components in $[-10^{-3}, 10^{-3}]$. For convenience, we set the EMCS evolution time sequence as $\boldsymbol{\tau} = [10^{-13}, 10^{-12}, 10^{-10}, 10^{-9}, 10^{-9}, 10^{-8}, 10^{-8}, 10^{-7}, 10^{-7}]~\rm{s}$. The hyper-parameters for the CSP method are set as $\text{tol}_{\text{rel}} = 10^{-3}$ and $\text{tol}_{\text{abs}} = 10^{-5}$. We then present the characteristic timescale distribution for each EMCS step. In Fig. \ref{fig:csp_time_preadtor_modulator}(B), it can be observed that the EMCS method effectively covers a broad range of timescales.

\subsection{\red{Range Estimation}}
\begin{figure}[htbp]
	\centering
        \includegraphics[width=0.9\textwidth]{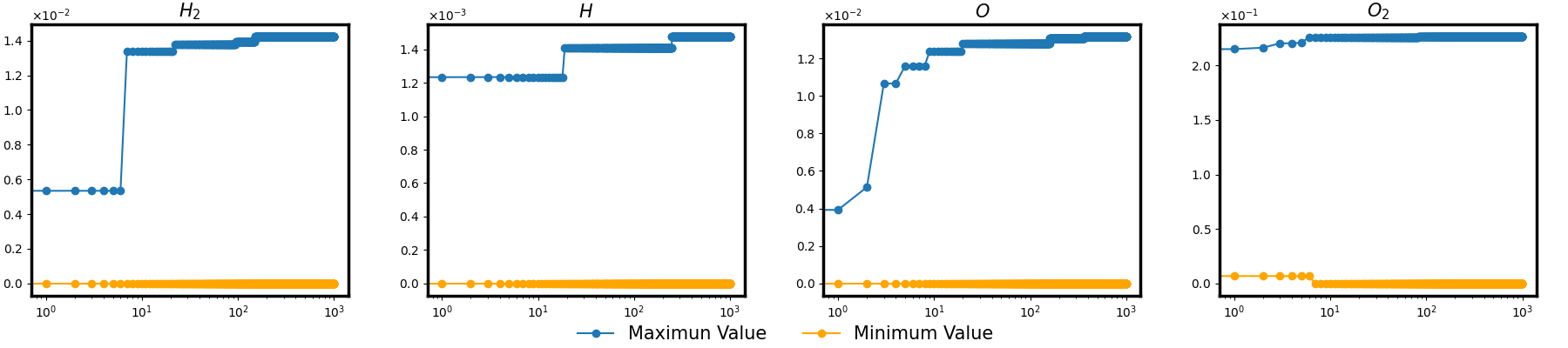}
        \caption{  \red{The estimated range of Step 1 in EMCS as the number of initial trajectories varies. The orange line represents the lower bound, and the blue line represents the upper bound.}  }
    \label{fig:max-min}
    \end{figure}

\red{ To demonstrate that a limited number of trajectories can provide a reasonable estimation of the data range, we visualized the dynamics of the maximum and minimum values for each dimension relative to the number of sampled trajectories. As shown in Fig.~\ref{fig:max-min}, the boundaries of the estimated range converge rapidly with only a small number of samples. The results reveal that a small number of data points based on manifold sampling could effectively estimate the data range.}

\subsection{\red{Ablation Study for DRM19 System}}
\label{sec:ablation}

\begin{figure}[htbp]
    	\centering
            \includegraphics[width=0.9\textwidth]{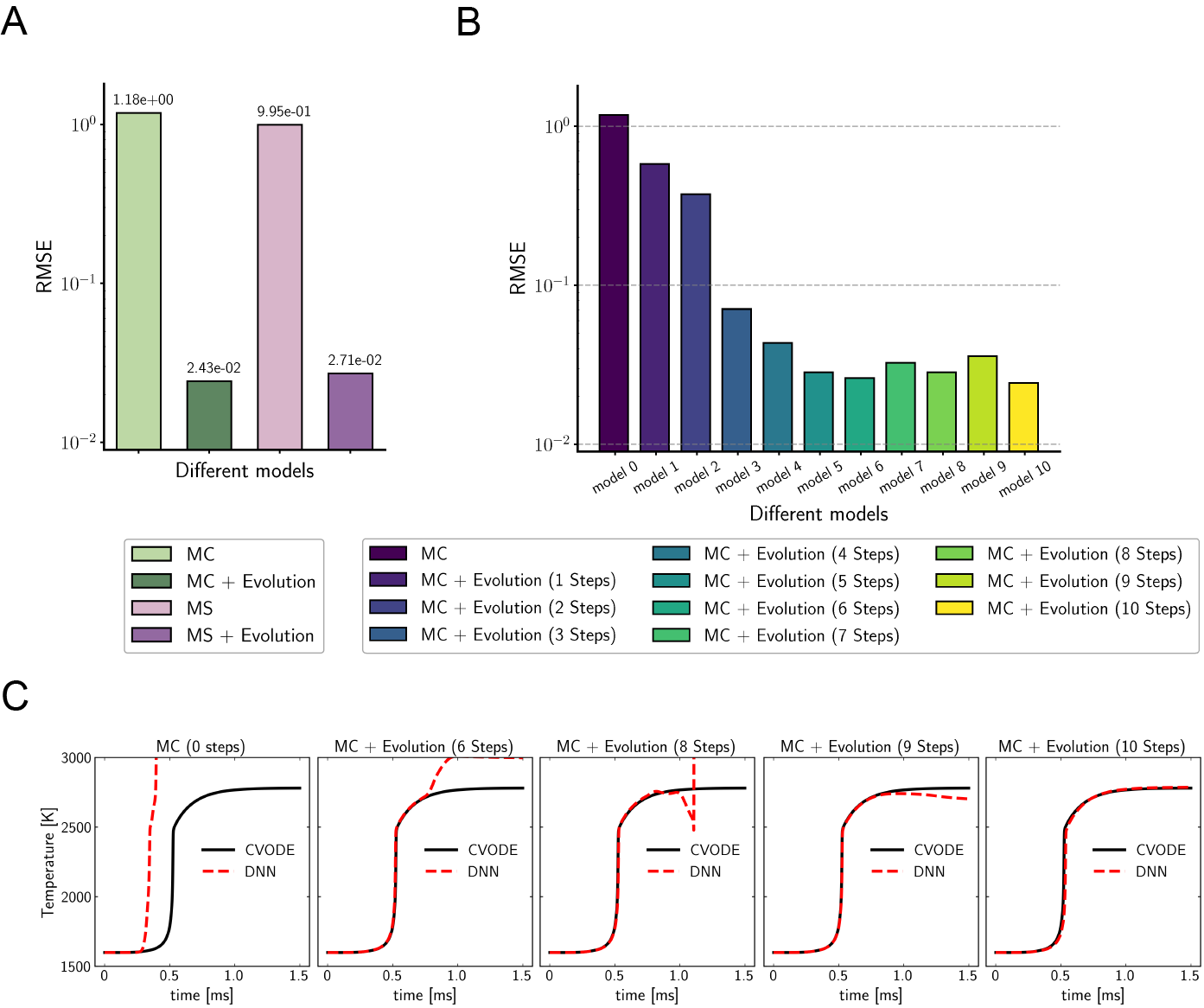}
    \caption{\red{ Ablation study of EMCS.  To ensure fairness, all four models in subplot (A) use the same data size of $1.1 \times 10^{7}$. The training set for the MC + Evolution model is generated by taking $10^{6}$ MC random samples and evolving them for 10 steps, while the training set for the MS + Evolution model is generated in a similar manner. The MC model in subplot (B) is consistent with that in subplot (A). The remaining models (labeled as model 1 through model 10) are generated by evolving $10^{6}$ MC initial samples for 1 to 10 steps, respectively.    }} 
    \label{fig:drm_ablation}
    \end{figure}

    \red{To quantify the contributions of the two EMCS components—Monte Carlo (MC) sampling and evolution-augmented generation—we conducted ablation experiments using the DRM19 chemical kinetics system (Section~\ref{sec:chemical}). We evaluated four models: MC sampling alone, MC with evolution augmentation (MC + Evolution), Multi-scale Sampling (MS)~\cite{zhang2022multi}, and MS with evolution augmentation (MS + Evolution), using the root mean square error (RMSE) of predicted temperatures as the metric. Fig.~\ref{fig:drm_ablation} (A) shows that baseline MC and MS models yield RMSEs of $\mathcal{O}(10^0)$, while their evolution-augmented counterparts reduce RMSE to $\mathcal{O}(10^{-2})$, highlighting the critical role of evolution augmentation. Fig.~\ref{fig:drm_ablation} (B, C) further demonstrates that increasing evolution steps (e.g., from two to three) significantly improve both point-wise and long-term temporal prediction accuracy, with diminishing returns beyond a few steps. These results confirm that evolution augmentation is the dominant contributor to model accuracy, as it captures the system’s multiscale dynamics, while MC sampling ensures broad phase space coverage.}

\subsection{\red{Baseline Comparisons with Advanced Sampling Techniques }}
\label{sec:baseline}

\begin{figure}[htbp]
    \centering
    \includegraphics[width=0.98\textwidth]{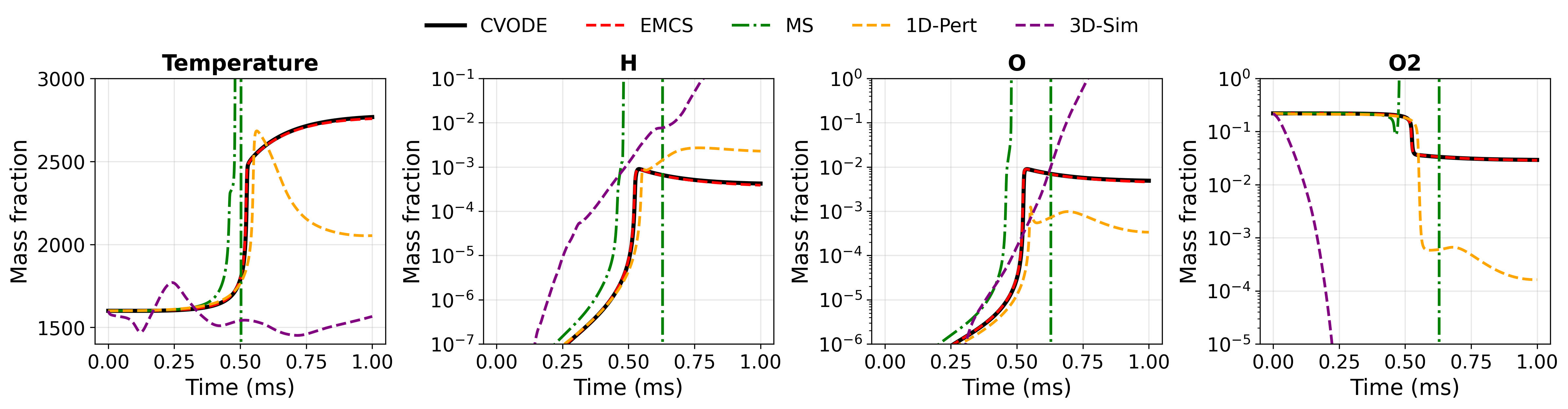}
    \caption{\red{Temporal prediction of models trained with different baseline sampling methods.}} 
    \label{fig:baseline_pred}
    \end{figure}

\begin{figure}[htbp]
    \centering
    \includegraphics[width=0.98\textwidth]{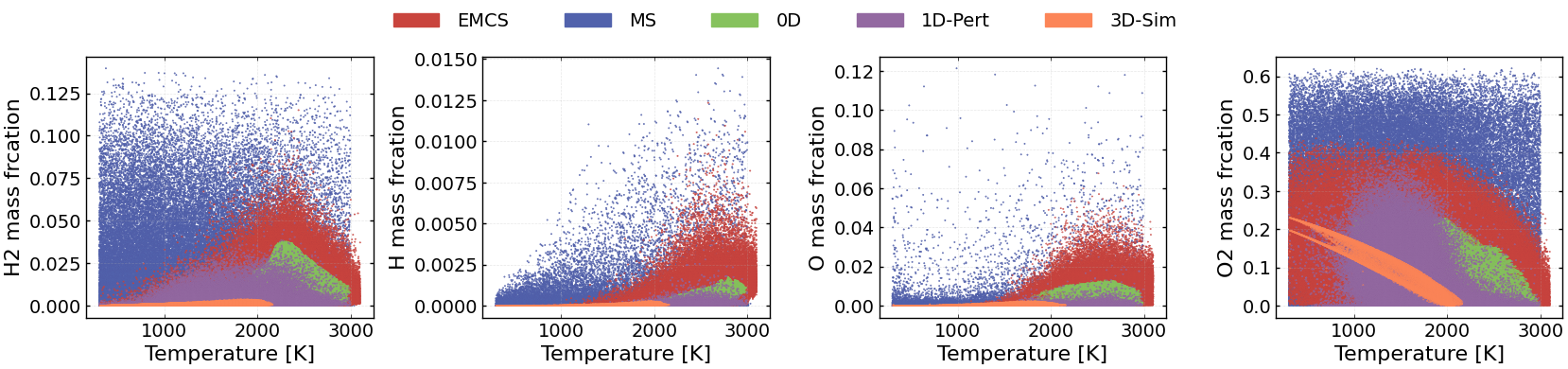}
    \caption{\red{The data distribution of different sampling methods.}} 
    \label{fig:baseline_dist}
    \end{figure}

  \red{
   To validate the advantages of the EMCS framework, we compared it against three baseline sampling techniques using the DRM19 chemical kinetics system (Section \ref{sec:chemical}): (1) direct sampling from 3D turbulent reactive flow simulations (3D-Sim)~\cite{weng2025CF, goswami2024CMAM}, (2) sampling from 1D steady-state reaction-diffusion systems with random perturbations (1D-Pert)~ \cite{ding2021CF, readshaw2023CF}, and (3) Multi-scale Sampling (MS) \cite{zhang2022multi}. Datasets were constructed as described in Section \ref{sec:chemical} and Appendix \ref{app:chemical}, and DNNs were trained to predict temporal evolution under identical initial conditions (1 atm, 1600 K, stoichiometric equivalence ratio). Fig.~\ref{fig:baseline_pred} shows that the EMCS-trained model demonstrates unique accuracy in its predictions, significantly outperforming 3D-Sim and MS model and maintaining stability unlike 1D-Pert, which diverges after 0.5 ms. Fig.~\ref{fig:baseline_dist} illustrates that 3D-Sim, 1D-Pert, and 0D datasets form distinct sub-manifolds with limited overlap, restricting generalization. EMCS, however, provides compact, configuration-independent coverage of the phase space, enhancing robustness across 0D, 1D, 2D, and 3D scenarios, as validated in Section \ref{sec:chemical} (Figs.~\ref{fig:tur_tri} and~\ref{fig:sandia}).}

\end{appendices}

\end{document}